\documentclass[a4paper,10pt]{article}
\usepackage{stmaryrd}
\usepackage{amsfonts}
\usepackage[latin1]{inputenc}
\usepackage{tikz-cd}

\usepackage{bbm}
\usepackage{amscd}
\usepackage{mathrsfs}
\usepackage{latexsym,amssymb,amsmath,amscd,amscd,amsthm,amsxtra}
\newtheorem{thm}{Theorem}[section]
\newtheorem{lem}[thm]{Lemma}
\newtheorem{cor}[thm]{Corollary}
\newtheorem{pro}[thm]{Proposition}

\newtheorem{rmk}[thm]{Remark}
\newtheorem{defi}[thm]{Definition}
\newtheorem{nota}[thm]{Notation}
\newtheorem{propdef}[thm]{Proposition-definition}

\newcommand {\emptycomment}[1]{}

\newcommand{\lon }{\,\rightarrow\,}
\newcommand{\be }{\begin{equation}}
\newcommand{\ee }{\end{equation}}

\newcommand{\pf}{\noindent{\bf Proof.}\ }

\newcommand{\g}{\frkg}
\newcommand{\h}{\frkh}

\newcommand{\R}{\mathbb R}
\newcommand{\LD}{\mathsf{Lie}\mathsf{Der}}
\newcommand{\Ld}{\mathsf{Lie2}\mathsf{Der}}
\newcommand{\Obg}{\mathsf{Ob^3_{(\omega_t,\varphi_t)}}}
\newcommand{\OBG}{\mathsf{Ob^2_{(\omega_t,\varphi_t)}}}
\newcommand{\der}{\mathsf{Der}}
\newcommand{\obe}{\mathsf{Ob^{\hat{\g}}_{(\varphi_\h,\varphi_\g)}}}
\newcommand{\rep}{\mathsf{Rep}}
\newcommand{\ob}{\mathsf{Ob}}

\newcommand{\huaB}{\mathcal{B}}

\newcommand{\huaV}{\mathcal{V}}

\newcommand{\huaC}{{\mathcal{C}}}

\newcommand{\huaH}{\mathcal{H}}

\newcommand{\huaZ}{\mathcal{Z}}

\newcommand{\frkg}{\mathfrak g}
\newcommand{\frkh}{\mathfrak h}

\newcommand{\frkk}{\mathfrak k}

\def\qed{\hfill ~\vrule height6pt width6pt depth0pt}

\newcommand{\br}[1]{   [ \cdot,    \cdot  ]_\frkg   }

\newcommand{\Id}{\rm{Id}}

\newcommand{\m}{\mathbbm{m}}

\newcommand{\dM}{\mathrm{d}}

\newcommand{\Hom}{\mathsf{Hom}}
\newcommand{\Der}{\mathsf{Der}}

\newcommand{\Ker}{\mathsf{ker}}

\newcommand{\gl}{\mathfrak {gl}}

\newcommand{\End}{\mathrm{End}}
\newcommand{\ad}{\mathsf{ad}}

\newcommand{\K}{\mathbb{K}}

\newcommand{\Lie}{\mathsf{Lie}}

\begin{document}
\title{
{Cohomologies of a Lie algebra with a derivation and applications}
\thanks
 {
Research supported by NSFC (11471139) and NSF of Jilin Province (20170101050JC).
 }
}
\author{Rong Tang, Yael Fr\'egier, Yunhe Sheng
 }

\date{}


\maketitle
\begin{abstract}
The main object of study of this paper is the notion of a $\LD$ pair, i.e. a Lie algebra with a derivation. We introduce the concept of a representation of a $\LD$ pair and study the corresponding cohomologies. We show that a $\LD$ pair is rigid if the second cohomology group is trivial, and a deformation of order $n$ is extensible if its obstruction class, which is defined to be an element is the third cohomology group, is trivial. We classify central extensions of $\LD$ pairs using the second cohomology group with the coefficient in the trivial representation. For a pair of derivations, we define its obstruction class and show that it is extensible if and only if the obstruction class is trivial. Finally, we classify $\Ld$ pairs using the third cohomology group of a $\LD$ pair.
\end{abstract}

\tableofcontents

\section{Introduction}

A classical approach to study a mathematical structure is to associate to it invariants. Among these, cohomology theories occupy a  central position as they enable for example to control deformations or extension problems. In particular, cohomology theories of various kinds of algebras have been developed with a great success (\cite{Ch-Ei,Gerstenhaber1, Har,Hor}). The deformation of algebraic structures began with the seminal work of Gerstenhaber~(\cite{Gerstenhaber2}) for associative algebras and followed by its extension to Lie algebras by Nijenhuis and Richardson~(\cite{Nijenhuis-Richardson}). After that, deformation of algebra morphisms and simultaneous deformations are widely studied (\cite{Fregier-Markl-Yau,Fregier,Fregier-Zambon1,Fregier-Zambon2,Gerstenhaber-Schack,Markl,Mandal,Yau}).

On the other hand, algebras are also useful via their derivations. For example one can mention their use in control theory (\cite{Ay-Ki-Tr,Ay-Ti}). One can also mention the higher derived bracket construction of Voronov in \cite{Vor} that produces out of some special derivation of a graded Lie algebra a homotopy Lie algebra. Moreover, algebras and their derivations have proven to be a very efficient tool to encode, via Koszul duality of operads of \cite{Gi-Ka}, many different types of structures as homological vector fields (formal geometry and Q-manifolds) (\cite{Sta, Vai}). They turn out to also play a fundamental role in the study of gauge theories in quantum field theory via the BV-formalism of \cite{Ba-Vi}. Indeed, this formalism relies on a Q-vector field $(S,\cdot)$ which is a derivation of the 1-shifted bracket $(\cdot,\cdot)$. In \cite{Loday}, Loday studied the operad of associative algebras with derivation.

 These facts  motivate us to construct in this paper a cohomology theory that controls, among other things, simultaneous deformations of a Lie algebra with a derivation.

 We begin by introducing in Section \ref{module} a convenient categorical framework by defining the category $\LD$ whose objects are Lie algebras with a derivation. The point is that this enables us to introduce modules and semi-direct products in this category as ``Beck-modules'', i.e. abelian group objects in the slice category. This is important since modules are a key ingredient in a cohomology theory.

 We then define in Section \ref{cohomology} a complex associated to a $\LD$ pair and  a module over it, and prove that it is indeed a complex. The rest of the paper is devoted to show that the usual interpretations of different cohomology groups are still valid in this framework. We first show in Section \ref{deformation} that the second cohomology group governs infinitesimal deformations modulo equivalences, and that obstructions to extension to higher order deformations are given by 3-cocycles. We then consider in Section \ref{extension} the classical problem of central extensions and its characterization in terms of second cohomology groups. In Section \ref{extension der} we study the extension problem of a pair of derivations. Finally we show in Section \ref{2-Lie} that third cohomology groups classify certain categorifications (skeletal 2-objects).

 There are many questions left open. The first one is to equip this complex with a compatible graded Lie algebra structure enabling to write the deformation equation as a Maurer-Cartan equation. It is probable that the dgLa can be obtained by a twist of a graded Lie algebra whose Maurer-Cartan elements are LieDer pairs. A by product of this approach would be to obtain at once the analogous theory for $L_{\infty}$ algebras with a (homotopy)-derivation. We expect to be able to obtain this more fundamental graded Lie algebra via the theory of operads, by showing that LieDer pairs are algebras over a suitable Koszul operad. Another line of research, suggested by the Beck-modules approach that we used to define our modules, could be to compare our complex with a complex that could be obtained by the Barr-Beck triple cohomology (\cite{Barr-Beck}). We intend to answer these questions in forth coming papers.

\vspace{2mm}
In this paper, we work over an algebraically closed field $\K$ of characteristic 0 and all the vector spaces are over $\K$.

\section{$\LD$ pairs and their modules}\label{module}

We need a suitable definition of module in the category of $\LD$ pairs. We follow Quillen's approach who characterized modules as monoid objects  in slice categories (\cite{Doubek-Markl-Zima}). This is why we start by introducing in Subsection \ref{LieDer} the category of interest for us. We then recall in Subsection \ref{Barr-Beck modules} the notion of slice category and monoid objects. We show that, given a Lie algebra $\g$, monoid objects in the slice category of Lie algebras over $\g$ are equivalent to $\g$-modules. This result is apparently well known to experts (\cite{Barr}), but we were not able to find its details in print. We then build on this result to obtain the $\LD$ version of module.

\subsection{The category of $\LD$ pairs} \label{LieDer}

Let $\g$ be a Lie algebra. A derivation of $\g$ is a linear map $\varphi:\g\lon\g$ which satisfies the Leibniz relation:
\begin{eqnarray}
\varphi[x,y]=[\varphi(x),y]+[x,\varphi(y)].
\end{eqnarray}
One denotes by $\Der(\g)$ the set of derivations of the Lie algebra $\g$.
\begin{defi}
A $\LD$ pair is a Lie algebra $\g$ with a derivation $\varphi\in\Der(\g)$. One denotes it by $(\g,\varphi)$. \end{defi}

It is called {\bf abelian} if the Lie bracket of $\g$ is trivial,
that is $[x,y]=0$ for all $x,y\in\g$.

\begin{nota}
We will use $[\cdot,\cdot]_\g$ instead of $[\cdot,\cdot]$  if precision is needed. The same goes for $\varphi_\g$ instead of $\varphi$.
\end{nota}

\begin{defi}
 Let $(\g,\varphi_\g)$ and $(\h,\varphi_\h)$ be $\LD$ pairs. A {\bf morphism} $f$ from $(\g,\varphi_\g)$ to $(\h,\varphi_\h)$ is a Lie algebra morphism $f:\g\lon \h$ such that
 \begin{eqnarray}
 f\circ\varphi_\g&=&\varphi_\h\circ f.
 \end{eqnarray}
\end{defi}
We denote by {\bf $\LD$} the category of $\LD$ pairs and their morphisms.

\subsection{Monoid objects in slice categories}\label{Barr-Beck modules}

\begin{defi}
For a category $\huaC$ and an object $A$ in $\huaC$. The {\bf slice category} $\huaC/A$ is the category whose
\begin{itemize}
\item objects $(B,\pi)$ are $\huaC$-morphisms $\pi:B\lon A,B\in\huaC$, and
\item morphisms $(B',\pi')\stackrel{f}{\lon}(B'',\pi'')$ are commutative diagrams of $\huaC$-morphisms:
\[\begin{CD}
B'@>f>>B''\\
  @V\pi' VV   @VV\pi''V  \\
A@=A.
\end{CD}\]

\end{itemize}
\end{defi}

\begin{defi}
 Let $\huaC$ be a category with finite products and a terminal object $T$. A {\bf monoid} object in $\huaC$ is an object $X\in \ob(\huaC)$ together with two morphisms $\mu:X\times X\lon X$ and $\eta:T\lon X$ such that following diagrams commute:
\begin{itemize}
\item the {\bf associativity} of $\mu:$

\[\begin{CD}\label{neut}
X\times X\times X@>\mu\times {\Id}_X>>X\times X\\
  @V{\Id}_X\times \mu VV   @VV\mu V  \\
X\times X@>>\mu >X,
\end{CD}\]
\item the {\bf neutrality} of $\eta:$
\[\begin{CD}
X\times X@>\mu>>           X        @<\mu<<                             X\times X        \\
  @A{\Id}_X\times \eta AA    @|                                   @AA\eta\times {\Id}_X A    \\
X\times T@<<({{\Id}_X},t_X) < X        @>>(t_X,{{\Id}_X})>                    T\times X,
\end{CD}\]
\end{itemize}
where $t_X:X\lon T$ is the unique map.

Let $\huaC_\m$ be the category whose objects are monoid objects $(X,\mu,\eta)$ in $\huaC$ as above and the hom-set $\Hom_{\huaC_\m}((X,\mu,\eta),(X',\mu',\eta'))$ is the set of all $f\in\Hom_{\huaC}(X,X')$ for which $\mu'\circ (f\times f)=f\circ \mu$ and $\eta'=f\circ\eta.$
\end{defi}



\subsection{The Lie case}

Denote by {\bf$\Lie$} the category of Lie algebras. Let $\g$ be a Lie algebra. We show in Subsection \ref{deci} how monoid objects in $\Lie/\g$ give rise to $\g$-modules. We then show in the following subsection that they form equivalent categories.

\subsubsection{Deciphering the definition of a monoid object in the Lie case}\label{deci}

 It is obvious that the terminal object of the slice category $\Lie/\g$ is $T=\g\stackrel{\Id}{\lon}\g$. Let $X=\frkk \stackrel{\kappa}{\lon}\g$ be a monoid object in $\Lie/\g$ with $\mu$ and $\eta$ as above. The information contained in $\eta$ gives this first result:

\begin{lem}\label{split}
There exists a vector space $V$ such that $$\frkk =s(\g)\oplus V.$$
\end{lem}

\pf It suffices to prove that $\kappa$ is a splitting epimorphism and take $V:= \Ker (\kappa)$. But remark that $\eta: T\lon X$ in the slice category actually means that there exists a commutative diagram

\begin{center}
  \begin{tikzpicture}[normal line/.style={->},font=\scriptsize]
\matrix (m) [matrix of math nodes, row sep=2em,
column sep=1.5em, text height=1.5ex, text depth=0.25ex]
{\g & &  \frkk\\
 & \g, & \\ };
\path[normal line]
(m-1-1) edge node[above] {$s$}(m-1-3)
edge node[below] {$\Id$} (m-2-2)
(m-1-3) edge node[below] {$\kappa$} (m-2-2);
\end{tikzpicture}
\end{center}
i.e. $\kappa\circ s=\Id_{}$, so we have $\frkk = s(\g)\oplus V$.  \qed\vspace{3mm}




The rest of this section aims to show that $V$ is an abelian Lie algebra and that it comes with an action of $\g$. But for that we need an expression of  $\mu$ in terms of $s$, the section of $\kappa$ that appeared in the previous proof.

First, notice that $X\times X$ is given by  $\frkk\times_{\g}\frkk \stackrel{\bar{\kappa}}{\lon}\g,$ where $\frkk\times_{\g}\frkk=\{(t,t')\in\frkk\oplus\frkk\mid \kappa(t)=\kappa(t')\},$ with bracket defined by $[(t,t'),(r,r')]:=([t,r],[t',r'])$ and $\bar{\kappa}(t,t'):=\kappa(t)$. In the slice category, $X\times X \stackrel{\mu}{\lon} X$ amounts to the commutative diagram

\begin{center}
  \begin{tikzpicture}[normal line/.style={->},font=\scriptsize]
\matrix (m) [matrix of math nodes, row sep=2em,
column sep=1.5em, text height=1.5ex, text depth=0.25ex]
{\frkk\times_{\g}\frkk & &  \frkk\\
 & \g. & \\ };
\path[normal line]
(m-1-1) edge node[above] {$M$}(m-1-3)
edge node[below] {$\bar{\kappa}$} (m-2-2)
(m-1-3) edge node[below] {$\kappa$} (m-2-2);
\end{tikzpicture}
\end{center}

\begin{lem}\label{Mul}With the above notations, we have
\begin{equation}
M(t,t')=t+t'-s(\kappa(t)).\label{Md}
\end{equation}
\end{lem}

\pf
In terms of $M$, the neutrality of $\eta$ amounts to the set of equations
\begin{eqnarray*}
M\circ (s\times \Id)\circ (\kappa,\Id)&=&\Id\\
M\circ (\Id\times s)\circ (\Id,\kappa)&=&\Id.
\end{eqnarray*}
We used here the fact that in the slice category, $t_X:X\lon T$ translates in the commutative diagram
\begin{center}
\begin{tikzpicture}[normal line/.style={->},font=\scriptsize]
\matrix (m) [matrix of math nodes, row sep=2em,
column sep=1.5em, text height=1.5ex, text depth=0.25ex]
{\frkk & & \g \\
 & \g & \\ };
\path[normal line]
(m-1-1) edge node[above] {$\kappa$}(m-1-3)
edge node[below] {$\kappa$} (m-2-2)
(m-1-3) edge node[below] {\quad$\Id$} (m-2-2);
\end{tikzpicture}
\end{center}

Thus, for any $t\in\frkk$ we have
\begin{eqnarray}\label{M}
M(s(\kappa(t)),t)=M(t,s(\kappa(t)))=t.
\end{eqnarray}
Let $(t,t')\in\frkk\times_{\g}\frkk$, since $\frkk\times_{\g}\frkk$ is a vector space, we have
\begin{eqnarray*}
(t,t')&=&\big(t-s(\kappa(t))+s(\kappa(t)),t'-s(\kappa(t'))+s(\kappa(t'))\big)\\
      &=&(s(\kappa(t)),s(\kappa(t')))+(t-s(\kappa(t)),0)+(0,t'-s(\kappa(t'))).
\end{eqnarray*}
Since $M$ is a morphism in $\Lie$, it is in particular linear, hence
\begin{eqnarray*}
M(t,t')=M(s(\kappa(t)),s(\kappa(t')))+M(t-s(\kappa(t)),0)+M(0,t'-s(\kappa(t'))).
\end{eqnarray*}
By \eqref{M}, we have
\begin{eqnarray*}
\nonumber M(t,t')&=&s(\kappa(t))+t-s(\kappa(t))+t'-s(\kappa(t'))\\
       &=&t+t'-s(\kappa(t)).
\end{eqnarray*} The proof is finished. \qed

\begin{rmk}
It can be seen from \eqref{Md} that a $\mu$ satisfying this equation is automatically associative. Therefore, a monoid object
in $\Lie/\g$ is uniquely determined by the neutral map $\eta$.
\end{rmk}

\begin{lem}
$V:=\Ker (\kappa)$ is an abelian sub Lie algebra of $\frkk$.
\end{lem}

\pf Since $\frkk\times_{\g}\frkk\stackrel{M}{\lon}\frkk$ is a map in $\Lie$,
\begin{eqnarray*}
M([(k,0),(0,k')])=[M(k,0),M(0,k')]=[k,k']
\end{eqnarray*}
by (\ref{Md}) if moreover $k,k'\in\Ker (\kappa)$.
On the left side, we have $$M([(k,0),(0,k')])=M([k,0],[0,k'])=M(0,0)=0.$$ Thus, we have $[k,k']=0$ for all $k,k'\in\Ker (\kappa)$.\qed\vspace{3mm}

We recall that a representation of a Lie algebra $\g$ on a vector space $V$ is a Lie algebra morphism $\rho:\g\lon\gl(V)$.

\begin{lem}\label{repres}
The expression
\begin{equation}\label{action}
\rho(x)(k)=[s(x),k]_{\frkk},\quad \forall x\in\g,~k\in V
\end{equation} defines a representation of $\g$ on $V$.
\end{lem}

\pf
Let us first check that $\rho$ is well defined, i.e. that for all $x\in\g,k\in V$, the right hand side of \eqref{action} is in $V$. It follows from
\begin{eqnarray*}
\kappa([s(x),k]_{\frkk})  =[x,\kappa(k)]_\g=[x,0]_\g=0,
\end{eqnarray*}
and $V=\Ker(\kappa)$.
We  now show that $\rho$ is a morphism. Recall that, by the Jacobi identity,  $\ad : \frkk\lon \gl(\frkk)$ is a Lie algebra morphism. One concludes by the remarks that $\rho=\ad\circ s$ and that $s$ is also a Lie algebra morphism.
\qed

\begin{rmk}\label{hint} Hadn't we known the definition of a representation, studying the properties of $\rho$ of Lemma \ref{repres} would have led us to the correct definition. This is what we will do in Section \ref{decild} to obtain the definition of a representation of a $\LD$ pair.
\end{rmk}

\subsubsection{Monoids in $\Lie/\g$ and $\g$-representations form equivalent categories}

We build two functors  \begin{tikzcd}
{(\Lie/\g)}_\m \arrow[bend left=35]{r}[name=Ker]{\Ker} & \g-\rep \arrow[bend left=35]{l}[name=ltimes]{\ltimes},
\end{tikzcd}and show that they induce equivalences of categories.

We start with the functor $\Ker:{(\Lie/\g)}_\m\lon \g$-$\rep$. Given $X=\frkk \stackrel{\kappa}{\lon}\g$
in ${(\Lie/\g)}_\m$, one defines $\Ker(X)$ to be $V:=\Ker (\kappa)$. The previous section insures that $\Ker(X)$ is in $\g$-$\rep$. We leave it as an exercise for the reader to verify that the restrictions of morphisms to the kernels of the maps to $\g$ induce maps of representations.

We now construct the functor $\ltimes : \g$-$\rep\lon {(\Lie/\g)}_\m$. Recall \cite{Humphreys} that one can associate to a representation $V$ of $\g$ the semi-direct product $\g\ltimes V$:
\begin{propdef}\label{semi}  Let $\rho:\g\lon\gl(V)$ be a representation of $\g$. The space $\g\oplus V$ becomes a Lie algebra with the bracket  $$[x+u,y+v]_{\rho}=[x,y]_{\g}+\rho(x)(v)-\rho(y)(u),\quad \forall x,y\in\g,~u,v\in V.$$ We denote this Lie algebra by $\g\ltimes V$.
\end{propdef}

 A first remark is that since the canonical projection $p :\g\ltimes V\lon \g$ is a morphism of Lie algebras, $X:=\g\ltimes V\stackrel{p}{\lon}\g$ is an object in $\Lie/\g$. It therefore suffices to equip $X$ with a suitable product $\mu$ and unit $\eta$ to complete the definition of the functor $\ltimes$.

 As already seen in the proof of Lemma \ref{split}, a map $\eta: T\lon X$ amounts to a splitting of the map $X=\frkk \stackrel{\kappa}{\lon}\g$. In the  case $X=\g\ltimes V\stackrel{p}{\lon}\g$, such a splitting, and hence our map $\eta$, is given by the canonical inclusion $i :\g \lon \g\ltimes V$.

  By Lemma \ref{Mul}, $\mu$ is determined by $M:(\g\ltimes V)\times_{\g}(\g\ltimes V)\lon\g\ltimes V$, where $(\g\ltimes V)\times_{\g}(\g\ltimes V)=\{(x+u,x+v)|x\in\g,u,v\in V\}$ and $M(x+u,x+v)=x+u+v$.

  Therefore, $\g\ltimes V\stackrel{p}{\lon}\g$, with these $\eta$ and $\mu$ is a monoid
object in $\Lie/\g$, the image $\ltimes(V)$ of $V$ by the functor $\ltimes$, which achieves the construction of this functor.

\emptycomment{
By Proposition \ref{tri}, we have a semi-direct product $\g\ltimes V$ $((\g\oplus V,[\cdot,\cdot]_{\rho});\varphi_\g+\varphi_V)$ of $((\g,[\cdot,\cdot]_\g);\varphi_\g))$ by $(\rho,V,\varphi_V)$. Let $G((\rho,V,\varphi_V))=\big(((\g\oplus V,[\cdot,\cdot]_{\rho});\varphi_\g+\varphi_V),p,\iota,\Ker (p)\big)$, where $p(x+u)=x$ and $\iota(x)=x$ for all $x\in\g$ and $u\in V$. From what has been discussed above, we have
$$X=\big(((\g\oplus V,[\cdot,\cdot]_{\rho});\varphi_\g+\varphi_V)\stackrel{p}{\lon}((\g,[\cdot,\cdot]_\g);\varphi_\g)\big)$$
 is a monoid object in the category $(\LD/((\g,[\cdot,\cdot]_\g);\varphi_\g))$. We can define a functor $G:\g$-$\rep$$\lon(\LD/((\g,[\cdot,\cdot]_\g);\varphi_\g))_{\m}$.
}

\begin{thm}\label{equi}
The functors  \begin{tikzcd}
{(\Lie/\g)}_\m \arrow[bend left=35]{r}[name=Ker]{\Ker} & \g-\rep \arrow[bend left=35]{l}[name=ltimes]{\ltimes}
\end{tikzcd} induce an equivalence of categories.
\end{thm}

\pf
It is straightforward to see that $\Ker\circ \ltimes=\Id$. On the other hand, for $\frkk \stackrel{\kappa}{\lon}\g\in(\Lie/\g)_\m$, we have
$$(\ltimes\circ \Ker)(\frkk \stackrel{\kappa}{\lon}\g)=(\g\ltimes \Ker (\kappa)\stackrel{p}{\lon}\g).$$
We define $\Psi:\frkk\lon \g\oplus \Ker (\kappa)$ by $\Psi(t)=\kappa(t)+t-(s\circ\kappa)(t)$. Moreover, $\Psi$ is a Lie algebra isomorphism from $\frkk$ to $\g\oplus \Ker (\kappa)$ and \begin{eqnarray*}
 p(\Psi(t))&=&p(\kappa(t)+t-(s\circ\kappa)(t))=\kappa(t),\\
 \Psi(s(x))&=&\kappa(s(x))+s(x)-(s\circ\kappa)(s(x))=x=i(x).
\end{eqnarray*}
Thus, $\Psi$ is an isomorphism from the monoid object $\frkk \stackrel{\kappa}{\lon}\g$ to the monoid object $\g\ltimes \Ker (\kappa)\stackrel{p}{\lon}\g$. We have defined a natural equivalence from the identify functor to $\ltimes\circ \Ker$.  \qed

\subsection{The $\LD$ case}

There is a forgetful functor $\LD\stackrel{F}{\lon}\Lie$ which consists in forgetting the derivations and that the maps intertwine the derivations.  Fix an object $(\g,\varphi_\g)$ in $\LD$. The functor $F$ induces a functor
$${(\LD/(\g,\varphi_\g))}_\m\stackrel{F}{\lon}{(\Lie/\g)}_\m$$ between  monoid objects in the slice categories.

The aim of this section is to lift the equivalence of the previous section at the level of  $\LD$ by completing the following diagram

\begin{center}

 \begin{tikzpicture}[normal line/.style={->},font=\scriptsize]
\matrix (m) [matrix of math nodes, row sep=2em,
column sep=1.5em, text height=1.5ex, text depth=0.25ex]
{{(\LD/(\g,\varphi_\g))}_\m & &  ?\\
& &\\
{(\Lie/\g)}_\m &  & \g-\rep.\\ };
\path[dotted](m-1-1) edge [bend left=35] (m-1-3)
(m-1-3) edge [bend left=35] (m-1-1)
edge  (m-3-3);
\path[->]
(m-1-1) edge node[left] {$F$} (m-3-1)
(m-3-1) edge [bend left=35] node[above] {$\Ker$} (m-3-3)
(m-3-3) edge [bend left=35] node[below] {$\ltimes$} (m-3-1);
\end{tikzpicture}
\end{center}
That is, by introducing the correct definition of module in the category $\LD$ and understanding its relationship with monoid objects in $\LD/(\g,\varphi_\g)$.

\subsubsection{Deciphering the definition of a monoid object in the $\LD$ case}\label{decild}

 According to Remark \ref{hint}, an analysis of monoid objects in $\LD/(\g,\varphi_\g)$ should lead us to the correct definition of module. We therefore consider an element $X=(\frkk,\varphi_{\frkk}) \stackrel{\kappa}{\lon}(\g,\varphi_\g)$ in ${(\LD/(\g,\varphi_\g))}_\m$ with $\mu$ and $\eta$ as above.

 The information contained in $\eta$ gives this first result:

\begin{lem}\label{splitphi}
One can decompose $\varphi_{\frkk}$ as
$$\varphi_{\frkk}=\varphi_{\frkk}|_{s(\g)}\oplus \varphi_{\frkk}|_{V}$$
where $s(\g)$ and $V$ are given by the decomposition $\frkk = s(\g)\oplus V$, see Lemma \ref{split}.
\end{lem}

\pf The morphism $\eta: T\lon X$ is given by the commutative diagram

\begin{center}
  \begin{tikzpicture}[normal line/.style={->},font=\scriptsize]
\matrix (m) [matrix of math nodes, row sep=2em,
column sep=1.5em, text height=1.5ex, text depth=0.25ex]
{(\g,\varphi_\g) & &  (\frkk,\varphi_{\frkk})\\
 & (\g,\varphi_\g). & \\ };
\path[normal line]
(m-1-1) edge node[above] {$s$}(m-1-3)
edge node[below] {$\Id$} (m-2-2)
(m-1-3) edge node[below] {$\kappa$} (m-2-2);
\end{tikzpicture}
\end{center}
Since $(\frkk,\varphi_{\frkk}) \stackrel{\kappa}{\lon}(\g,\varphi_\g)$ is in $\LD$, that is, $\kappa\circ\varphi_{\frkk}=\varphi_\g\circ \kappa$. For all $v\in \Ker(\kappa)=V$, we have
$$
\kappa(\varphi_{\frkk}(v))=\varphi_\g(\kappa(v))=0.
$$
Thus, we deduce that $\varphi_{\frkk}(V)\subset V$. Similarly, $(\g,\varphi_\g) \stackrel{s}{\lon}(\frkk,\varphi_{\frkk})$ is in $\LD$, for all $x\in\g$, we have
$$
\varphi_{\frkk}(s(x))=s(\varphi_\g(x)).
$$
Thus, we deduce that $\varphi_{\frkk}(s(\g))\subset s(\g)$. The proof is finished.
\qed\vspace{3mm}

We now interpret in terms of $\varphi_{\g}$ and $\varphi_{\frkk}|_{V}$ the fact that $\varphi_{\frkk}$ is a derivation.

\begin{lem}\label{leib}
The following is satisfied
$$\varphi_{\frkk}|_{V}[s(x),v]=[s(\varphi_{\frkg}(x)),v]+[s(x),\varphi_{\frkk}|_{V}(v)].$$
\end{lem}

\pf
Apply the Leibniz rule for $\varphi_{\frkk}=\varphi_{\frkk}|_{s(\g)}\oplus \varphi_{\frkk}|_{V}$ to the element $[s(x),v]$ of $\frkk = s(\g)\oplus V$ and $\varphi_{\frkk}\circ s=s\circ\varphi_{\frkg}$.
\qed\vspace{3mm}

We can now, according to Remark \ref{hint} and the previous two lemmas, give the following :

\begin{defi} \label{represLD}
A {\bf representation} of a $\LD$ pair $(\g,\varphi_\g)$ on a vector space $V$ with respect to $\varphi_V\in\mathfrak{gl}(V)$ is a Lie algebra morphism $\rho:\mathfrak{g}\lon \mathfrak{gl}(V)$ such that for all $x\in\mathfrak{g}$, the following equality is satisfied:
\begin{eqnarray}
\label{rep1}\varphi_V\circ \rho(x)&=&\rho(\varphi_\g(x))+\rho(x)\circ\varphi_V.
\end{eqnarray}
\end{defi}
We denote a representation  by $(\rho,V,\varphi_V)$. For all $x\in\mathfrak{g}$, we define $\ad_{x}:\mathfrak{g}\lon \mathfrak{g}$ by
\begin{eqnarray}
\ad_{x}(y)=[x,y]_{\mathfrak{g}},\quad\forall y \in \mathfrak{g}.
\end{eqnarray}
Then $(\ad,\g,\varphi_\g)$ is a representation of the $\LD$ pair $(\g,\varphi_\g)$
on $\g$ with respect to $\varphi_\g$, which is called the {\bf adjoint representation}.

A representation $(\rho,V,\varphi_V)$ of a $\LD$ pair $(\g,\varphi_\g)$ is said to be {\bf trivial} if $\rho=0$.

A similar study of morphisms of monoid objects leads to the following
\begin{defi}
Let $(\rho,V,\varphi_V)$ and $(\rho',V',\varphi_{V'})$ be two representations of the $\LD$ pair $(\g,\varphi_\g)$. A {\bf morphism} from $(\rho,V,\varphi_V)$ to $(\rho',V',\varphi_{V'})$ is a morphism of Lie algebra representations $f:V\lon V'$  such that
\begin{eqnarray}
f\circ\varphi_V=\varphi_{V'}\circ f.
\end{eqnarray}
\end{defi}

\begin{nota}
Let $(\g,\varphi_\g)$ be a $\LD$ pair. We denote by $(\g,\varphi_\g)$-$\rep$ the category of the representations of the $\LD$ pair $(\g,\varphi_\g)$ and their morphisms.
\end{nota}

To sum-up the discussion and as an immediate corollary of Definition \ref{represLD} and Lemmas \ref{repres}, \ref{splitphi} and \ref{leib}, one has

\begin{cor}\label{Cor}
Given a monoid object $X=(\frkk,\varphi_{\frkk}) \stackrel{\kappa}{\lon}(\g,\varphi_\g)$ in $\LD/(\g,\varphi_\g)$, the expression $$\rho(x)(k)=[s(x),k]_{\frkk}$$ defines a representation of $(\g,\varphi_\g)$ on $V=\ker(\kappa)$ with respect to $\varphi_{\frkk}|_{V}\in\mathfrak{gl}(V)$.
\end{cor}

\subsubsection{Monoids in $\LD/(\g,\varphi_\g)$ and $(\g,\varphi_\g)$-representations form equivalent categories}

In the previous section we have partially completed the diagram, lifting in Corollary \ref{Cor} the functor $\Ker$ to the $\LD$ level

\begin{center}

 \begin{tikzpicture}[normal line/.style={->},font=\scriptsize]
\matrix (m) [matrix of math nodes, row sep=2em,
column sep=1.5em, text height=1.5ex, text depth=0.25ex]
{{(\LD/(\g,\varphi_\g))}_\m & &  (\g,\varphi_\g)-\rep\\
& &\\
{(\Lie/\g)}_\m &  & \g-\rep\\ };
\path[dotted]
(m-1-3) edge [bend left=35] (m-1-1);
\path[->]
(m-1-3)edge node[right] {$F$} (m-3-3)
(m-1-1) edge node[left] {$F$} (m-3-1)edge [bend left=35] node[above] {$\Ker$} (m-1-3)
(m-3-1) edge [bend left=35] node[above] {$\Ker$} (m-3-3)
(m-3-3) edge [bend left=35] node[below] {$\ltimes$} (m-3-1);
\end{tikzpicture}

\end{center}

We now focus on the task of showing that the functor $\ltimes$ can also be lifted, i.e. we want to show that it is compatible with derivations.

\begin{pro}\label{tri}
Given a representation $(\rho,V,\varphi_V)$ of a $\LD$ pair $(\g,\varphi_\g)$, define $\varphi_\g+\varphi_V:\g\oplus V\lon\g\oplus V$ by
$$(\varphi_\g+\varphi_V)(x+u)=\varphi_\g(x)+\varphi_V(u).$$
Then $(\g\oplus V,\varphi_\g+\varphi_V)$, with the Lie structure of Proposition \ref{semi},  is a $\LD$ pair which we call the semi-direct product of the $\LD$ pair $(\g,\varphi_\g)$ by the representation $(\rho,V,\varphi_V)$ and denote it by $\g\ltimes_{\LD} V$.
\end{pro}

\pf It suffices to show that $\varphi_\g+\varphi_V$ is a derivation. On one hand, we have
\begin{eqnarray*}
(\varphi_\g+\varphi_V)[x+u,y+v]_{\rho}&=&(\varphi_\g+\varphi_V)([x,y]_{\g}+\rho(x)(v)-\rho(y)(u))\\
                                      &=&\varphi_\g([x,y]_{\g})+(\varphi_V\circ\rho(x))(v)-(\varphi_V\circ\rho(y))(u).
\end{eqnarray*}
On the other hand, we have
\begin{eqnarray*}
&&[(\varphi_\g+\varphi_V)(x+u),y+v]_{\rho}+[x+u,(\varphi_\g+\varphi_V)(y+v)]_{\rho}\\
=&&[\varphi_\g(x)+\varphi_V(u),y+v]_{\rho}+[x+u,\varphi_\g(y)+\varphi_V(v)]_{\rho}\\
=&&[\varphi_\g(x),y]_\g+\rho(\varphi_\g(x))(v)-(\rho(y)\circ\varphi_V)(u)+[x,\varphi_\g(y)]_\g+(\rho(x)\circ\varphi_V)(v)-\rho(\varphi_\g(y))(u).
\end{eqnarray*}
 One concludes by \eqref{rep1} and the fact that $\varphi_\g$ is a derivation. \qed

\emptycomment{
By Proposition \ref{tri}, we have a semi-direct product $\g\ltimes V$ $((\g\oplus V,[\cdot,\cdot]_{\rho});\varphi_\g+\varphi_V)$ of $((\g,[\cdot,\cdot]_\g);\varphi_\g))$ by $(\rho,V,\varphi_V)$. Let $G((\rho,V,\varphi_V))=\big(((\g\oplus V,[\cdot,\cdot]_{\rho});\varphi_\g+\varphi_V),p,\iota,\Ker (p)\big)$, where $p(x+u)=x$ and $\iota(x)=x$ for all $x\in\g$ and $u\in V$. From what has been discussed above, we have
$$X=\big(((\g\oplus V,[\cdot,\cdot]_{\rho});\varphi_\g+\varphi_V)\stackrel{p}{\lon}((\g,[\cdot,\cdot]_\g);\varphi_\g)\big)$$
 is a monoid object in the category $(\LD/((\g,[\cdot,\cdot]_\g);\varphi_\g))$. We can define a functor $G:\g$-$\rep$$\lon(\LD/((\g,[\cdot,\cdot]_\g);\varphi_\g))_{\m}$.
}

\begin{thm}
The functors  \begin{tikzcd}
{(\LD/(\g,\varphi_\g))}_\m \arrow[bend left=35]{r}[name=Ker]{\Ker} &   (\g,\varphi_\g)-\rep
 \arrow[bend left=35]{l}[name=G]{\ltimes_{\LD}}
\end{tikzcd} induce an equivalence of categories.
\end{thm}

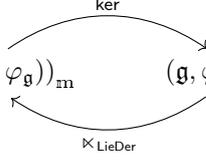
\pf By Theorem \ref{equi},  we only need to prove that
$$
\Psi\circ\varphi_{\frkk}=(\varphi_{\g}\oplus \varphi_{\frkk}|_{V})\circ\Psi.
$$
Since $s$ and $\kappa$ are $\LD$ pairs morphisms. For all $t\in\frkk$, we have
\begin{eqnarray*}
(\varphi_{\g}\oplus \varphi_{\frkk}|_{V})(\Psi(t))&=&(\varphi_{\g}\oplus \varphi_{\frkk}|_{V})\big(\kappa(t)+t-(s\circ\kappa)(t)\big)\\
&=&\varphi_{\g}(\kappa(t))+\varphi_{\frkk}(t-(s\circ\kappa)(t))\\
&=&\kappa(\varphi_{\frkk}(t))+\varphi_{\frkk}(t)-s(\varphi_{\g}(\kappa(t)))\\
&=&\kappa(\varphi_{\frkk}(t))+\varphi_{\frkk}(t)-s(\kappa(\varphi_{\frkk}(t)))\\
&=&(\Psi\circ\varphi_{\frkk})(t).
\end{eqnarray*}
The proof is finished. \qed

\section{Cohomologies of $\LD$ pairs} \label{cohomology}

Let $(\rho,V)$ be a representation of a Lie algebra $\g$. The Chevalley-Eilenberg cohomology of the Lie algebra $\g$ with the coefficient in $V$ is the cohomology of the cochain complex $C^n(\g;V)=\Hom(\wedge^{n}\mathfrak{g},V)$  with the
coboundary operator $\dM : C^n(\g;V )\lon C^{n+1}(\g;V)$ defined by
\begin{eqnarray*}
(\dM f_n)(x_1,\cdots,x_{n+1})&=&\sum_{i=1}^{n+1}(-1)^{i+1}\rho(x_{i})(f(x_1,\cdots,\widehat{x_{i}},\cdots,x_{n+1}))\\
                          &&+\sum_{i<j}(-1)^{i+j}f([x_i,x_j],x_1,\cdots,\widehat{x_i},\cdots,\widehat{x_j},\cdots,x_{n+1}).
\end{eqnarray*}
Denoted the set of closed $n$-cochains by $\huaZ^n(\g;V)$ and the set of exact $n$-cochains by $\huaB^n(\g;V)$. We denote by $\huaH^n(\g;V)=\huaZ^n(\g;V)/\huaB^n(\g;V)$ the corresponding cohomology group.

Let $(\rho,V,\varphi_V)$  be a representation of a $\LD$ pair $(\g,\varphi_\g)$.
We define the set of $\LD$ pair $0$-cochains to be $0$, and define
the set of $\LD$ pair $1$-cochains to be $C^1_{\LD}(\g;V)=\Hom(\g,V)$. For $n\geq 2$, we define the set of $\LD$ pair $n$-cochains by
\begin{eqnarray*}
C^n_{\LD}(\g;V)\stackrel{\triangle}{=}C^{n}(\g;V)\times C^{n-1}(\g;V).
\end{eqnarray*}

For $n\geq1$, we define an operator $\delta:C^n(\g;V)\lon C^n(\g;V)$ by
\begin{eqnarray*}
\delta f_n=\sum_{i=1}^{n}f_n\circ({\Id\otimes\cdots\otimes\varphi_\g\otimes\cdots\otimes Id})-
 \varphi_V \circ f_n.
\end{eqnarray*}

\emptycomment{
Define $\partial:C^0_{\LD}(\g;V)\longrightarrow C^1_{\LD}(\g;V)$ by
\begin{equation}
  \partial (u)=\dM u,\quad\forall u\in V,
\end{equation}
and define $\partial:C^1_{\LD}(\g;V)\longrightarrow C^2_{\LD}(\g;V)$ by
\begin{equation}
  \partial (\theta,u)=(\dM \theta,\delta \theta),\quad\forall \theta\in \Hom(\g,V).
  \end{equation}
  }
Define $\partial:C^1_{\LD}(\g;V)\longrightarrow C^2_{\LD}(\g;V)$ by
\begin{equation}
  \partial f_1=(\dM f_1,(-1)^1\delta f_1),\quad\forall f_1\in \Hom(\g,V).
  \end{equation}
Then for $n\geq2,$ we define $\partial:C^n_{\LD}(\g;V)\lon C^{n+1}_{\LD}(\g;V)$ by
\begin{eqnarray}
\partial(f_{n},g_{n-1})=\Big({\dM} f_{n},{\dM} g_{n-1}+(-1)^{n}\delta f_n \Big),\quad \forall f_n\in C^n(\g;V),~g_{n-1}\in C^{n-1}(\g;V).
\end{eqnarray}

The following lemma gives the relation between the operator $\dM$ and the operator $\delta$, which plays important role in the proof of that $\partial$ is a coboundary operator. We omit the proof which is straightforward tedious computations.

\begin{lem}\label{coboundary oprator}
The map $\dM$ and $\delta$ are commutative with each other, i.e. $\dM\circ\delta=\delta\circ\dM$.
\end{lem}

\begin{thm}
The map $\partial$ is a coboundary operator, i.e. $\partial\circ\partial$=$0$.
\end{thm}
\pf For $n\geq1,$ since $\dM\circ\delta=\delta\circ\dM$, we have
\begin{eqnarray*}
(\partial\circ\partial)(f_{n},g_{n-1})&=&\partial(\dM f_{n},\dM g_{n-1}+(-1)^{n}\delta f_{n})\\
                                            &=&\Big(\dM(\dM f_{n}),\dM(\dM g_{n-1}+(-1)^{n}\delta f_{n})+(-1)^{n+1}\delta(\dM f_{n})\Big)\\
                                            &=&\Big(0,(-1)^{n}(\dM\circ\delta-\delta\circ\dM)(f_{n})\Big)\\
                                            &=&0.
\end{eqnarray*}
 Thus,   the map $\partial$ is a coboundary operator. \qed\vspace{3mm}

Associated to the representation $(\rho,V,\varphi_V)$, we obtain a complex $(C^*_{\LD}(\g;V),\partial)$. Denoted the set of closed $n$-cochains by $\huaZ_{\LD}^n(\g;V)$ and the set of exact $n$-cochains by $\huaB_{\LD}^n(\g;V)$. We define the corresponding cohomology group by
$$\huaH_{\LD}^n(\g;V)=\huaZ_{\LD}^n(\g;V)/\huaB_{\LD}^n(\g;V).$$


\begin{pro}
Let $(\rho,V,\varphi_V)$ be a representation of a $\LD$ pair $(\g,\varphi_\g)$. Then we have
\begin{eqnarray*}
\huaH_{\LD}^1(\g;V)&=&\{f|f\in\huaZ^1(\g;V),~ f\circ\varphi_\g=\varphi_V\circ f\},
\end{eqnarray*}
 \end{pro}
\pf   For any $f\in C^1_{\LD}(\g;V)$, we have
\begin{eqnarray*}
\partial f=(\dM f,-\delta f).
\end{eqnarray*}
Therefore, $f$ is closed if and only if $f\in \huaZ^1(\g;V)$ and $f\circ\varphi_\g=\varphi_V\circ f$. The conclusion follows from the fact that there is no exact $1$-cochain. \qed


\section{Deformations of a $\LD$ pair }\label{deformation}
In this section, we study formal deformations and deformations of order $n$ of a $\LD$ pair.

\subsection{Formal deformations of a $\LD$ pair }

In this subsection, we study $1$-parameter formal deformations of a $\LD$ pair. We show that if the second cohomology group of a $\LD$ pair with the coefficient in the adjoint representation is trivial, then the $\LD$ pair is rigid.

Let $(\g,\varphi)$ be a $\LD$ pair. In the sequel, we will also denote the Lie bracket $[\cdot,\cdot]$ by $\omega$. Consider a $t$-parametrized family of linear operations
\begin{eqnarray*}
\omega_t&=&\sum_{i\geq0}\omega_it^i,\quad \omega_i\in C^2(\g;\g),\\
\varphi_t&=&\sum_{i\geq0}\varphi_i t^i,\quad \varphi_i\in C^1(\g;\g).
\end{eqnarray*}

\begin{defi}
 If all $(\omega_t,\varphi_t)$ endow the $\mathbb{K}[[t]]$-module $\g[[t]]$\footnote{The notation $\g[[t]]$ means the vector space of formal power series in $t$ with coefficients in $\g$, that is, for all $x_t\in\g[[t]]$, we have $x_t=x_0+x_1t+x_2t^2+\cdots,$  for $x_0,x_1,x_2,\cdots\in\g.$} the $\LD$-pair structure with $(\omega_0,\varphi_0)=(\omega,\varphi)$, we say that $(\omega_t,\varphi_t)$
is a {\bf  $1$-parameter formal deformation} of the $\LD$ pair $(\g,\varphi)$.
\end{defi}
A pair $(\omega_t,\varphi_t)$, as given above, is a $1$-parameter formal deformation of the $\LD$ pair $(\g,\varphi)$ if and only if for all $x,y,z\in \g$, the following equalities   hold:
\begin{eqnarray}
\label{deformation1}&&\omega_t(\omega_t(x,y),z)+\omega_t(\omega_t(y,z),x)+\omega_t(\omega_t(z,x),y)=0,\\
\label{deformation2}&&\varphi_t(\omega_t(x,y))-\omega_t(\varphi_t(x),y)-\omega_t(x,\varphi_t(y))=0.
\end{eqnarray}
 Expanding the equations in \eqref{deformation1}, \eqref{deformation2} and collecting coefficients of $t^n$, we see that  \eqref{deformation1} and \eqref{deformation2} are equivalent to the system of equations
\begin{eqnarray}
\label{deformation3}&&\sum\limits_{i+j=n\atop i,j\geq0}\Big(\omega_i(\omega_j(x,y),z)+\omega_i(\omega_j(y,z),x)+\omega_i(\omega_j(z,x),y)\Big)=0,\\
\label{deformation4}&&\sum\limits_{i+j=n\atop i,j\geq0}\Big(\varphi_i(\omega_j(x,y))-\omega_j(\varphi_i(x),y)-\omega_j(x,\varphi_i(y))\Big)=0.
\end{eqnarray}

\begin{rmk}\label{rmk}
For $n=0$, condition \eqref{deformation3} is equivalent to the usual Jacobi identity of $\omega$, and \eqref{deformation4} is equivalent to the fact that $\varphi$ is a derivation.
 \end{rmk}

 \begin{pro}\label{pro:cocycle}
 Let $(\omega_t,\varphi_t)$
be a $1$-parameter formal deformation of the $\LD$ pair $(\g,\varphi)$.
Then $(\omega_1,\varphi_1)\in C^2_{\LD}(\g;\g)$   is a $2$-cocycle of the $\LD$ pair  $(\g,\varphi)$ with the coefficient in the adjoint representation.
\end{pro}
\pf
 For $n=1$, \eqref{deformation3} is equivalent to $\dM \omega_1=0$, and \eqref{deformation4} is equivalent to $\dM\varphi_1+\delta \omega_1=0$. Thus for $n=1$, \eqref{deformation3} and \eqref{deformation4} are equivalent to  that $(\omega_1,\varphi_1)$ is a $2$-cocycle.  \qed

\begin{defi}
The $2$-cocycle $(\omega_1,\varphi_1)$ is called the {\bf infinitesimal} of the $1$-parameter formal deformation $(\omega_t,\varphi_t)$ of the $\LD$ pair $(\g,\varphi)$. 
\end{defi}


\begin{defi}
Let $(\bar{\omega}_t,\bar{\varphi}_t)$ and $(\omega_t,\varphi_t)$ be $1$-parameter formal deformations of a $\LD$ pair $(\g,\varphi)$. A
{\bf formal isomorphism} from $(\bar{\omega}_t,\bar{\varphi}_t)$ to $(\omega_t,\varphi_t)$ is a power series $\Phi_t=\sum_{i\geq0}\phi_it^i:\g[[t]]\lon\g[[t]]$, where  $\phi_i\in C^1(\g;\g)$ with $\phi_0=\Id$, such that
\begin{eqnarray}
&&\Phi_t\circ\bar{\omega}_t= \omega_t\circ(\Phi_t\times\Phi_t),\\
&&\Phi_t\circ\bar{\varphi}_t=\varphi_t\circ\Phi_t.
\end{eqnarray}


Two $1$-parameter formal deformations $(\bar{\omega}_t,\bar{\varphi}_t)$ and $(\omega_t,\varphi_t)$ are said to be {\bf equivalent} if  there exists a formal isomorphism $\Phi_t:(\bar{\omega}_t,\bar{\varphi}_t)\lon(\omega_t,\varphi_t)$.
\end{defi}

\begin{thm}
The infinitesimals of two equivalent $1$-parameter formal deformations of a $\LD$ pair $(\g,\varphi)$ are in the same cohomology class.
\end{thm}
\pf Let $\Phi_t:(\bar{\omega}_t,\bar{\varphi}_t)\lon(\omega_t,\varphi_t)$ be a formal isomorphism. Then for all $x,y\in \g$ we have
\begin{eqnarray*}
\Phi_t\circ\bar{\omega}_t(x,y)&=& \omega_t\circ(\Phi_t\times\Phi_t)(x,y),\\
\Phi_t\circ\bar{\varphi}_t(x)&=& \varphi_t\circ\Phi_t (x).
\end{eqnarray*}
 Expanding the above identities and comparing coefficients of $t$, we have
\begin{eqnarray*}
\bar{\omega}_1(x,y)&=&\omega_1(x,y)+\omega(\phi_1(x),y)+\omega(x,\phi_1(y))-\phi_1(\omega(x,y)),\\
\bar{\varphi}_1(x)&=&\varphi_1(x)+\varphi(\phi_1(x))-\phi_1(\varphi(x)).
\end{eqnarray*}
Thus, we have $(\bar{\omega}_1,\bar{\varphi}_1)=(\omega_1,\varphi_1)+\partial(\phi_1)$, which implies that $[(\bar{\omega}_1,\bar{\varphi}_1)]=[(\omega_1,\varphi_1)]\in\huaH^2_{\LD}(\g;\g)$. The proof is finished. \qed 


\begin{defi}
A $1$-parameter formal deformation $(\omega_t,\varphi_t)$ of a $\LD$ pair $(\g,\varphi)$ is said to be {\bf trivial} if it is equivalent to $(\omega,\varphi)$, i.e. there exists  $\Phi_t=\sum_{i\geq0}\phi_it^i:\g[[t]]\lon\g[[t]]$, where  $\phi_i\in C^1(\g;\g)$ with $\phi_0=\Id$, such that
\begin{eqnarray}
&&\Phi_t\circ\bar{\omega}_t= \omega\circ(\Phi_t\times\Phi_t),\\
&&\Phi_t \circ\bar{\varphi}_t=\varphi\circ\Phi_t.
\end{eqnarray}
 \end{defi}

\begin{defi}
A $\LD$ pair $(\g,\varphi)$ is said to be {\bf rigid} if  every $1$-parameter formal deformation of $(\g,\varphi)$ is trivial. 
\end{defi}

\begin{thm}
If $\huaH^2_{\LD}(\g;\g)=0$, then the $\LD$ pair $(\g,\varphi)$ is rigid.
\end{thm}
\pf Let $(\omega_t,\varphi_t)$ be a $1$-parameter formal deformation of the $\LD$ pair $(\g,\varphi)$. By Proposition \ref{pro:cocycle}, we have $(\omega_1,\varphi_1)\in\huaZ^2_{\LD}(\g;\g)$. By $\huaH^2_{\LD}(\g;\g)=0,$ there exists a 1-cochain $\phi_1\in  C^1_{\LD}(\g;\g)$ such that
\begin{eqnarray}
\label{rigid}(\omega_1,\varphi_1)=-\partial(\phi_1).
\end{eqnarray}
 Then setting $\Phi_t={\Id}+\phi_1 t$, we have a deformation $(\bar{\omega}_t,\bar{\varphi}_t)$, where
\begin{eqnarray*}
\bar{\omega}_t(x,y)&=&\big(\Phi_t^{-1}\circ \omega_t\circ(\Phi_t\times\Phi_t)\big)(x,y),\\
\bar{\varphi}_t(x)&=&\big(\Phi_t^{-1}\circ\varphi_t\circ\Phi_t\big)(x).
\end{eqnarray*}
Thus, $(\bar{\omega}_t,\bar{\varphi}_t)$ is equivalent to $(\omega_t,\varphi_t)$. Moreover, we have
\begin{eqnarray*}
\bar{\omega}_t(x,y)&=&({\Id}-\phi_1t+\phi_1^2t^{2}+\cdots+(-1)^i\phi_1^it^{i}+\cdots)(\omega_t(x+\phi_1(x)t,y+\phi_1(y)t)),\\
\bar{\varphi}_t(x)&=&({\Id}-\phi_1t+\phi_1^2t^{2}+\cdots+(-1)^i\phi_1^it^{i}+\cdots)(\varphi_t(x+\phi_1(x)t)).
\end{eqnarray*}
Thus, we have
\begin{eqnarray*}
\bar{\omega}_t(x,y)&=&\omega(x,y)+(\omega_1(x,y)+\omega(x,\phi_1(y))+\omega(\phi_1(x),y)-\phi_1(\omega(x,y)))t+\bar{\omega}_{2}(x,y)t^{2}+\cdots,\\
\bar{\varphi}_t(x)&=&\varphi(x)+(\varphi(\phi_1(x))+\varphi_1(x)-\phi_1(\varphi(x)))t+\bar{\varphi}_{2}(x)t^{2}+\cdots.
\end{eqnarray*}
By \eqref{rigid}, we have
\begin{eqnarray*}
\bar{\omega}_t(x,y)&=&\omega(x,y)+\bar{\omega}_{2}(x,y)t^{2}+\cdots,\\
\bar{\varphi}_t(x)&=&\varphi(x)+\bar{\varphi}_{2}(x)t^{2}+\cdots.
\end{eqnarray*}
Then by repeating the argument, we can show that $(\omega_t,\varphi_t)$ is equivalent to $(\omega,\varphi)$. \qed
\subsection{Deformations of order $n$ of a $\LD$ pair}
In this subsection, we introduce a cohomology class associated to any deformation of order $n$ of a $\LD$ pair. We show that a deformation of order $n$ of a $\LD$ pair is extensible if and only if this cohomology class is trivial. Thus we call this cohomology class the obstruction class of a deformation of order $n$ being extensible.

\begin{defi}
A deformation of order $n$ of a $\LD$ pair $(\g,\varphi)$ is a pair $(\omega_t,\varphi_t)$ such that $\omega_t=\sum_{i=0}^n\omega_it^i$ and $\varphi_t=\sum_{i=0}^n\varphi_i t^i$ endow the $\mathbb{K}[[t]]/(t^{n+1})$-module $\g[[t]]/(t^{n+1})$ the $\LD$ pair structure with $(\omega_0,\varphi_0)=(\omega,\varphi)$.
\end{defi}

\begin{defi}
Let $(\omega_t,\varphi_t)$ be a deformation of order $n$ of a $\LD$ pair $(\g,\varphi)$. If there exists a $2$-cochain $(\omega_{n+1},\varphi_{n+1})\in C^2_{\LD}(\g;\g)$, such that the pair $(\widetilde{\omega}_t,\widetilde{\varphi}_t)$ with $$\widetilde{\omega}_t=\omega_t+\omega_{n+1}t^{n+1},\quad \widetilde{\varphi}_t=\varphi_t+\varphi_{n+1}t^{n+1}$$ is a deformation of order $n+1$ of the $\LD$ pair $(\g,\varphi)$, then we say that $(\omega_t,\varphi_t)$ is {\bf extensible}.
\end{defi}

Let $(\omega_t,\varphi_t)$ be a deformation of order $n$ of a $\LD$ pair $(\g,\varphi)$. Define $(\Obg,\OBG)\in C^3_{\LD}(\g;\g)$ by
\begin{eqnarray}
\label{eq:obs31}\Obg(x,y,z)&=&\sum\limits_{i+j=n+1\atop i,j>0}\big(\omega_i(\omega_{j}(x,y),z)+\omega_i(\omega_{j}(y,z),x)+\omega_i(\omega_{j}(z,x),y)\big),\\
\label{eq:obs32}\OBG(x,y)&=&\sum\limits_{i+j=n+1\atop i,j>0}\big(\varphi_i(\omega_{j}(x,y))-\omega_{j}(\varphi_i(x),y)-\omega_{j}(x,\varphi_i(y))\big).
\end{eqnarray}

In the sequel, we show that $(\Obg,\OBG)$ is a 3-cocycle. We need some preparations.
Let $\g$ be a Lie algebra. The Nijenhuis-Richardson bracket $[\cdot,\cdot]$  on
 the graded vector space $ C^*(\mathfrak{g};\mathfrak{g})=\oplus_{k=0}^{+\infty}C^{k+1}(\mathfrak{g};\mathfrak{g})$ is given by \cite{Nijenhuis-Richardson}:
\begin{eqnarray}
[P,Q]=P\circ Q-(-1)^{pq}Q\circ P,\,\,\,\,\forall P\in C^{p+1}(\mathfrak{g};\mathfrak{g}),
Q\in C^{q+1}(\mathfrak{g};\mathfrak{g}),
\end{eqnarray}
where  $P\circ Q\in C^{p+q+1}(\mathfrak{g};\mathfrak{g})$ is defined by
\begin{eqnarray}
P\circ Q(x_{1},\cdots,x_{p+q+1})=\sum_{\sigma\in unsh(q+1,p)}(-1)^{\sigma}
P(Q(x_{\sigma(1)},\cdots,x_{\sigma(q+1)}),
x_{\sigma(q+2)},\cdots,x_{\sigma(p+q+1)}).
\end{eqnarray}
Furthermore, $(C^*(\mathfrak{g};\mathfrak{g}),[\cdot,\cdot],\bar{\partial}=[\omega,\cdot])$ is a differential graded Lie algebra. The Chevalley-Eilenberg coboundary operator $\dM$ of the Lie algebra $\g$ with the coefficient in the adjoint representation can be reformulated as follows:
\begin{eqnarray}
\dM f=(-1)^{k-1}\bar{\partial}f=(-1)^{k-1}[\omega,f],\quad \forall f\in C^k(\g;\g).
\end{eqnarray}

\begin{pro}
Let $(\omega_t,\varphi_t)$ be a deformation of order $n$ of a $\LD$ pair $(\g,\varphi)$. The $3$-cochain $(\Obg,\OBG)$ defined by \eqref{eq:obs31} and \eqref{eq:obs32} is a $3$-cocycle of the $\LD$ pair $(\g,\varphi)$ with the coefficient in the adjoint representation.
\end{pro}

\pf We use the Nijenhuis-Richardson bracket to write $\Obg$ and $\OBG$  as follows:
\begin{eqnarray}
\Obg=\frac{1}{2}\sum\limits_{i+j=n+1\atop i,j>0}[\omega_i,\omega_j],\quad
\OBG=\sum\limits_{i+j=n+1\atop i,j>0}[\varphi_i,\omega_j].
\end{eqnarray}
Since $(\omega_t,\varphi_t)$ is a deformation of order $n$ of the $\LD$ pair $(\g,\varphi)$, for all $0\leq i\leq n$, we have
\begin{eqnarray}
\label{deformation3.1}&&\sum\limits_{k+l=i\atop k,l\geq0}\Big(\omega_k(\omega_l(x,y),z)+\omega_k(\omega_l(y,z),x)+\omega_k(\omega_l(z,x),y)\Big)=0,\\
\label{deformation4.1}&&\sum\limits_{k+l=i\atop k,l\geq0}\Big(\varphi_k(\omega_l(x,y))-\omega_l(\varphi_k(x),y)-\omega_l(x,\varphi_k(y))\Big)=0.
\end{eqnarray}
Thus, the equations \eqref{deformation3.1} and \eqref{deformation4.1} are equivalent to
\begin{eqnarray}
\label{Lie}&&\frac{1}{2}\sum\limits_{k+l=i\atop k,l>0}[\omega_k,\omega_l]=-[\omega,\omega_i],\\
\label{deri}&&\sum\limits_{k+l=i\atop k,l>0}[\varphi_k,\omega_l]=-[\varphi,\omega_i]+[\omega,\varphi_i].
\end{eqnarray}
Then we have
\begin{eqnarray*}
\dM \Obg&=&(-1)^2[\omega,\Obg]\\
       &=&\frac{1}{2}\sum\limits_{i+j=n+1\atop i,j>0}[\omega,[\omega_i,\omega_j]]\\
       &=&\frac{1}{2}\sum\limits_{i+j=n+1\atop i,j>0}\big([[\omega,\omega_i],\omega_j]-[\omega_i,[\omega,\omega_j]]\big)\\
       &\stackrel{\eqref{Lie}}{=}&-\frac{1}{4}\sum\limits_{i'+i''+j=n+1\atop i',i'',j>0}[[\omega_{i'},\omega_{i''}],\omega_j]+\frac{1}{4}
       \sum\limits_{i+j'+j''=n+1\atop i,j',j''>0}[\omega_i,[\omega_{j'},\omega_{j''}]]\\
       &=&\frac{1}{4}\sum\limits_{i'+i''+j=n+1\atop i',i'',j>0}[\omega_j,[\omega_{i'},\omega_{i''}]]+\frac{1}{4}
       \sum\limits_{i+j'+j''=n+1\atop i,j',j''>0}[\omega_i,[\omega_{j'},\omega_{j''}]]\\
       &=&\frac{1}{2}\sum\limits_{i'+i''+j=n+1\atop i',i'',j>0}[\omega_j,[\omega_{i'},\omega_{i''}]]\\
       &=&0.
\end{eqnarray*}
Moreover, for all $f\in C^k(\g;\g)$ we have
\begin{eqnarray}
\delta f=-[\varphi, f].
\end{eqnarray}
Thus, we have
\begin{eqnarray*}
&&\dM \OBG+(-1)^3\delta\Obg\\&=&-[\omega, \OBG]+[\varphi,\Obg]\\
                                     &=&-\sum\limits_{i+j=n+1\atop i,j>0}[\omega,[\varphi_i,\omega_j]]+\frac{1}{2}\sum\limits_{i+j=n+1\atop i,j>0}[\varphi,[\omega_i,\omega_j]]\\
                                     \nonumber&=&-\sum\limits_{i+j=n+1\atop i,j>0}\big([[\omega,\varphi_i],\omega_j]+[\varphi_i,[\omega,\omega_j]]\big)+\frac{1}{2}\sum\limits_{i+j=n+1\atop i,j>0}
   \big([[\varphi,\omega_i],\omega_j]+[\omega_i,[\varphi,\omega_j]]\big)\\
                                     &=& -\sum\limits_{i+j=n+1\atop i,j>0}\big([[\omega,\varphi_i],\omega_j]+[\varphi_i,[\omega,\omega_j]]\big)+\sum\limits_{i+j=n+1\atop i,j>0}
   [[\varphi,\omega_i],\omega_j]\\
    &\stackrel{\eqref{Lie}}{=}&-\sum\limits_{i+j=n+1\atop i,j>0}[[\omega,\varphi_i],\omega_j]+\frac{1}{2}
       \sum\limits_{i+j'+j''=n+1\atop i,j',j''>0}[\varphi_i,[\omega_{j'},\omega_{j''}]]
       +\sum\limits_{i+j=n+1\atop i,j>0}
   [[\varphi,\omega_i],\omega_j]\\
   &\stackrel{\eqref{deri}}{=}&-\sum\limits_{i'+i''+j=n+1\atop i',i'',j>0}[[\varphi_{i'},\omega_{i''}],\omega_j]-\sum\limits_{i+j=n+1\atop i,j>0}[[\varphi,\omega_i],\omega_j]
   +\frac{1}{2}\sum\limits_{i+j'+j''=n+1\atop i,j',j''>0}[[\varphi_i,\omega_{j'}],\omega_{j''}]\\
   &&+\frac{1}{2}\sum\limits_{i+j'+j''=n+1\atop i,j',j''>0}[\omega_{j'},[\varphi_i,\omega_{j''}]]
    +\sum\limits_{i+j=n+1\atop i,j>0}[[\varphi,\omega_i],\omega_j]\\
    &=&-\sum\limits_{i'+i''+j=n+1\atop i',i'',j>0}[[\varphi_{i'},\omega_{i''}],\omega_j]+\sum\limits_{i+j'+j''=n+1\atop i,j',j''>0}[[\varphi_i,\omega_{j'}],\omega_{j''}]\\
    &=&0.
\end{eqnarray*}
Therefore, we have
$$\partial(\Obg,\OBG)=(\dM \Obg,\dM \OBG+(-1)^3\delta\Obg)=0.$$
  The proof is finished. \qed

  \begin{defi}
  Let $(\omega_t,\varphi_t)$ be a deformation of order $n$ of a $\LD$ pair $(\g,\varphi)$.  The cohomology class $[(\Obg,\OBG)]\in\huaH^3_\LD(\g;\g)$ is called the {\bf obstruction class} of  $(\omega_t,\varphi_t)$ being extensible.
  \end{defi}

\begin{thm}
Let $(\omega_t,\varphi_t)$ be a deformation of order $n$ of a $\LD$ pair $(\g,\varphi)$. Then $(\omega_t,\varphi_t)$ is extensible if and only if the obstruction class $[(\Obg,\OBG)]$ is trivial.
\end{thm}
\pf Suppose that a deformation $(\omega_t,\varphi_t)$ of order $n$ of the $\LD$ pair $(\g,\varphi)$ extends to a deformation of order $n+1$. Then \eqref{deformation3.1} and \eqref{deformation4.1} hold for $i=n+1$. Thus, we have
\begin{eqnarray}
\Obg=\dM \omega_{n+1},\quad
\OBG=\dM\varphi_{n+1}+\delta \omega_{n+1},
\end{eqnarray}
which implies that
\begin{eqnarray}
(\Obg,\OBG)=\partial(\omega_{n+1},\varphi_{n+1}).
\end{eqnarray}
Thus, the obstruction class $[(\Obg,\OBG)]$ is trivial.

 Conversely, if the obstruction class $[(\Obg,\OBG)]$ is trivial, suppose that
\begin{eqnarray}
(\Obg,\OBG)=\partial(\omega_{n+1},\varphi_{n+1}),
\end{eqnarray}
for some 2-cochain $(\omega_{n+1},\varphi_{n+1})\in C^2_{\LD}(\g;\g)$. Set
\begin{eqnarray}
(\widetilde{\omega}_t,\widetilde{\varphi}_t)=(\omega_t+\omega_{n+1}t^{n+1},\varphi_t+\varphi_{n+1}t^{n+1}).
\end{eqnarray}
Then $(\widetilde{\omega}_t,\widetilde{\varphi}_t)$ satisfies \eqref{deformation3.1}-\eqref{deformation4.1} for $0\leq i\leq n+1$, so $(\widetilde{\omega}_t,\widetilde{\varphi}_t)$ is a deformation of order $n+1$, which implies that  $(\omega_t,\varphi_t)$ is extensible. \qed

\begin{cor}
If $\huaH^3_{\LD}(\g;\g)=0$, then every $2$-cocycle in $C^2_{\LD}(\g;\g)$ is the infinitesimal of some $1$-parameter formal deformation of the $\LD$ pair $(\g,\varphi)$.
\end{cor}

\section{Central extensions of a $\LD$ pair}\label{extension}
In this section, we study central extensions of a $\LD$ pair and show that  central extensions of a $\LD$ pair $(\g,\varphi_\g)$ are controlled by the second cohomology of $(\g,\varphi_\g)$ with the coefficient in the trivial representation.

\begin{defi}
Let $(\h,\varphi_\h)$ be an abelian $\LD$ pair and $(\g,\varphi_\g)$ a $\LD$ pair. An exact
sequence of $\LD$ pair morphisms
\[\begin{CD}
0@>>>\mathfrak{h}@>\iota>>\hat{\mathfrak{g}}@>p>>{\mathfrak{g}}             @>>>0\\
@.    @V\varphi_\h VV   @V\varphi_{\hat{\mathfrak{g}}}VV  @V\varphi_\g VV    @.\\
0@>>>\mathfrak{h}@>\iota>>\hat{\mathfrak{g}}@>p>>{\mathfrak{g}}             @>>>0
\end{CD}\]
is called a {\bf central extension} of $(\g,\varphi_\g)$ by $(\h,\varphi_\h)$, if $[\h,\hat{\g}]_{\hat{\g}}=0$, that is $[h,x]_{\hat{\g}}=0$ for all $h\in\h$ and $x\in\hat{\g}$. Here we identify $\h$ with the corresponding
subalgebra of $\hat{\g}$. Therefore, we have $\varphi_{\hat{\mathfrak{g}}}|_{\h}=\varphi_\h.$
\end{defi}

\begin{defi}\label{defi:iso}
Let $(\hat{\g}_1,\varphi_{\hat{\g}_1})$ and $(\hat{\g}_2,\varphi_{\hat{\g}_2})$ be two central extensions of $(\g,\varphi_\g)$ by $(\h,\varphi_\h)$. They are said to be {\bf isomorphic} if there exists a $\LD$ pair morphism $\zeta:(\hat{\g}_1,\varphi_{\hat{\g}_1})\lon (\hat{\g}_2,\varphi_{\hat{\g}_2})$ such that we have the following commutative diagram:
\label{iso}\[\begin{CD}
0@>>>(\mathfrak{h},\varphi_\h)@>\iota_{1}>>(\hat{\g}_1,\varphi_{\hat{\g}_1})@>p_{1}>>{(\mathfrak{g},\varphi_\g)}             @>>>0\\
@.    @|                       @V\zeta VV                     @|                       @.\\
0@>>>(\mathfrak{h},\varphi_\h)@>\iota_{2}>>(\hat{\g}_2,\varphi_{\hat{\g}_2})@>p_{2}>>{(\mathfrak{g},\varphi_\g)}             @>>>0
.\end{CD}\]
\end{defi}

A {\bf section} of a central extension $(\hat{\g},\varphi_{\hat{\g}})$ of $(\g,\varphi_\g)$ by $(\h,\varphi_\h)$ is a linear map $s:\g\longrightarrow \hat{\g}$ such that $p\circ s=\Id$.

Let $(\hat{\g},\varphi_{\hat{\g}})$ be a central extension of a $\LD$ pair $(\g,\varphi_\g)$ by an abelian $\LD$ pair $(\h,\varphi_\h)$ and $s:\frkg\lon \hat{\frkg}$ a section.
Define linear maps $\psi:\mathfrak{g}\wedge\mathfrak{g}\lon \mathfrak{h}$ and $\chi:\g\lon\h$ respectively by
\begin{eqnarray}
\label{do}\psi(x,y)&=&[s(x),s(y)]_{\hat{\mathfrak{g}}}-s[x,y]_{\mathfrak{g}}, \,\,\,\,\forall x,y \in \mathfrak{g},\\
\label{dr}\chi(x)&=&\varphi_{\hat{\g}}(s(x))-s(\varphi_\g(x)), \,\,\,\,\forall x\in \mathfrak{g}.
\end{eqnarray}

Obviously, $\hat{\g}$ is isomorphic to $\g\oplus\h$ as vector spaces. Transfer the $\LD$ pair structure on $\hat{\mathfrak{g}}$ to that on $\mathfrak{g}\oplus \mathfrak{h}$, we obtain a $\LD$ pair $(\mathfrak{g}\oplus \mathfrak{h},\varphi_{\chi})$, where the Lie bracket $[\cdot,\cdot]_{\psi}$ and $\varphi_{\chi}$ are given by
\begin{eqnarray}
\label{dbr}[x+h,y+l]_{\psi}&=&[x,y]_{\mathfrak{g}}+\psi(x,y), \,\,\,\,\forall x,y \in \mathfrak{g},\,\,h,l\in\h,\\
 \label{dmo}\varphi_{\chi}(x+h)&=&\varphi_{\mathfrak{g}}(x)+\chi(x)+\varphi_{\mathfrak{h}}(h), \,\,\,\,\forall x\in \mathfrak{g},\,\,h\in\h.
\end{eqnarray}
The following proposition gives the conditions on $\psi$ and $\chi$ such that $(\g\oplus\h,\varphi_{\chi})$ is a $\LD$ pair.
\begin{pro}\label{lieder}
With the above notations, $(\g\oplus\h,\varphi_{\chi})$ is a $\LD$ pair if and only if $(\psi,\chi)$ is a $2$-cocycle of  the $\LD$ pair $(\g,\varphi_\g)$ with the coefficient in the trivial representation $(\rho=0,\frkh,\varphi_\frkh)$, i.e. $(\psi,\chi)$ satisfy the following equalities:
\begin{eqnarray}
\label{p1}&&\psi([x,y]_\g,z)+\psi([y,z]_\g,x)+\psi([z,x]_\g,y)=0,\\
\label{p2}&&\chi([x,y]_\g)+\varphi_\h(\psi(x,y))-\psi(\varphi_\g(x),y)-\psi(x,\varphi_\g(y))=0.
\end{eqnarray}
\end{pro}
\pf If $(\g\oplus\h,\varphi_{\chi})$ is a $\LD$ pair, by $$[[x+h,y+l]_{\psi},z+t]_{\psi}+c.p.=0,$$  we deduce that \eqref{p1} holds. By
$$\varphi_{\chi}([x+h,y+l]_{\psi})=[\varphi_{\chi}(x+h),y+l]_{\psi}+[x+h,\varphi_{\chi}(y+l)]_{\psi},$$
we deduce that \eqref{p2} holds.

Conversely, if \eqref{p1} and \eqref{p2} hold, it is straightforward to see that $(\g\oplus\h,\varphi_{\chi})$ is a $\LD$ pair. The proof is finished.\qed

\begin{thm}
Let $(\h,\varphi_\h)$ be an abelian $\LD$ pair and $(\g,\varphi_\g)$ a $\LD$ pair. Then central extensions of $(\g,\varphi_\g)$ by $(\h,\varphi_\h)$ are classified by the second cohomology group $\huaH^2_{\LD}(\g;\h)$ of  the $\LD$ pair $(\g,\varphi_\g)$ with the coefficient in the trivial representation $(\rho=0,\frkh,\varphi_\frkh)$.
\end{thm}
\pf Let $(\hat{\g},\varphi_{\hat{\g}})$ be a central extension of $(\g,\varphi_\g)$ by $(\h,\varphi_\h)$. By choosing a section $s:\mathfrak{g}\lon \hat{\mathfrak{g}}$, we obtain a $2$-cocycle $(\psi,\chi)$.  Now we show that the cohomological class of $(\psi,\chi)$ does not depend on the choice of sections. In fact, let $s_{1}$ and $s_{2}$ be two different sections. Define $\phi:\mathfrak{g}\lon \mathfrak{h}$ by $\phi(x)=s_1(x)-s_2(x)$. Then we have
\begin{eqnarray*}
\psi_{1}(x,y)&=&[s_{1}(x),s_{1}(y)]_{\hat{\mathfrak{g}}}-s_{1}[x,y]_{\mathfrak{g}}\\
&=&[s_{2}(x)+\phi(x),s_{2}(y)+\phi(y)]_{\hat{\mathfrak{g}}}-s_{2}[x,y]_{\mathfrak{g}}-\phi([x,y]_{\mathfrak{g}})\\
&=&[s_{2}(x),s_{2}(y)]_{\hat{\mathfrak{g}}}-s_{2}[x,y]_{\mathfrak{g}}-\phi([x,y]_{\mathfrak{g}})\\
&=&\psi_{2}(x,y)-\phi([x,y]_{\mathfrak{g}}),
\end{eqnarray*}
and
\begin{eqnarray*}
\chi_1(x)&=&\varphi_{\hat{\g}}(s_1x)-s_1(\varphi_\g(x))\\
         &=&\varphi_{\hat{\g}}(s_{2}(x)+\phi(x))-(s_{2}+\phi)(\varphi_\g(x))\\
         &=&\varphi_{\hat{\g}}(s_{2}(x))-s_{2}(\varphi_\g(x))+\varphi_\h(\phi(x))-\phi(\varphi_\g(x))\\
         &=&\chi_2(x)+\varphi_\h(\phi(x))-\phi(\varphi_\g(x)).
\end{eqnarray*}
Thus, we obtain $(\psi_{1},\chi_1)=(\psi_{2},\chi_2)+\partial(\phi)$. Therefore, $(\psi_{1},\chi_1)$ and $(\psi_{2},\chi_2)$ are in the same cohomological class.

Now we go on to prove that isomorphic central extensions give rise to the same element in $\huaH^2_{\LD}(\g;\h)$. Assume that $(\hat{\g}_1,\varphi_{\hat{\g}_1})$ and $(\hat{\g}_2,\varphi_{\hat{\g}_2})$ are two isomorphic central extensions of $(\g,\varphi_\g)$ by $(\h,\varphi_\h)$, and $\zeta:(\hat{\g}_1,\varphi_{\hat{\g}_1})\lon (\hat{\g}_2,\varphi_{\hat{\g}_2})$ is a $\LD$ pair morphism such that we have the commutative diagram in Definition \ref{iso}. Assume that $s_{1}:\mathfrak{g}\lon \hat{\mathfrak{g}}_{1}$ is a section of $\hat{\g}_1$. By $p_{2}\circ\zeta=p_{1}$, we have
$$p_{2}\circ(\zeta\circ s_{1})=p_{1}\circ s_{1}=\Id.$$
Thus, we obtain that $\zeta\circ s_{1}$ is   a section of $\hat{\g}_2$. Define $s_2=\zeta\circ s_{1}$. Since $\zeta$ is a morphism of $\LD$ pair and $\zeta\mid_{\mathfrak{h}}=\Id$, we have
\begin{eqnarray*}
\psi_2(x,y)&=&[s_{2}(x),s_{2}(y)]_{\hat{\mathfrak{g}}_{2}}-s_{2}[x,y]_{\mathfrak{g}}=
[(\zeta\circ s_{1})(x),(\zeta\circ s_{1})(y)]_{\hat{\mathfrak{g}}_{2}}-(\zeta\circ s_{1})[x,y]_{\mathfrak{g}}\\
&=&\zeta([s_{1}(x),s_{1}(y)]_{\hat{\mathfrak{g}}_{1}}-s_{1}[x,y]_{\mathfrak{g}})\\
&=&[s_{1}(x),s_{1}(y)]_{\hat{\mathfrak{g}}_{1}}-s_{1}[x,y]_{\mathfrak{g}}\\
&=&\psi_{1}(x,y),
\end{eqnarray*}
and
\begin{eqnarray*}
\chi_2(x)&=&\varphi_{\hat{\g}_2}(s_{2}x)-s_{2}(\varphi_\g(x))=\varphi_{\hat{\g}_2}((\zeta\circ s_{1})(x))-(\zeta\circ s_{1})(\varphi_\g(x))\\
                                                 &=&\zeta(\varphi_{\hat{\g}_1}(s_{1}x)-s_{1}(\varphi_\g(x)))\\
                                                 &=&\varphi_{\hat{\g}_1}(s_{1}x)-s_{1}(\varphi_\g(x))\\
                                                 &=&\chi_1(x).
\end{eqnarray*}
 Thus, the isomorphic central extensions give rise to the same element in $\huaH^2_{\LD}(\g;\h)$.

Conversely, given two 2-cocycles $(\psi_{1},\chi_1)$ and $(\psi_{2},\chi_2)$, we can construct two central extensions  $(\g\oplus\h,\varphi_{\chi_1})$ and $(\g\oplus\h,\varphi_{\chi_2})$, as in $\eqref{dbr}$ and $\eqref{dmo}$. If they represent the same cohomological class, i.e. there exists $\phi:\mathfrak{g}\lon \mathfrak{h}$, such that $(\psi_{1},\chi_1)=(\psi_{2},\chi_2)+\partial(\phi)$, we define $\zeta:\mathfrak{g}\oplus \mathfrak{h}\lon \mathfrak{g}\oplus \mathfrak{h}$ by
$$\zeta(x+h)=x+\phi(x)+h.$$
Then we can deduce that $\zeta$ is   an isomorphism between central extensions.
We omit details. This finishes the proof.\qed

\section{Extensions of a pair of derivations}\label{extension der}

In \cite{Bardakov-Singh}, the authors study extensions of a pair of automorphisms of Lie algebras. Since derivations are infinitesimals of automorphisms, we are interested in extensions of a pair of derivations. In this section, associated to a central extension $\hat{\g}$ of a Lie algebra $\g$ by an abelian Lie algebra $\h$ and a pair of derivations $(\varphi_\h,\varphi_\g)\in\Der(\h)\times\Der(\g)$, we define a cohomology class $[\obe]\in \huaH^2(\g;\h)$. We show that   $(\varphi_\h,\varphi_\g)$ is extensible if and only if the cohomology class $[\obe]$ is trivial. Thus we call $[\obe]$ the obstruction class of  $(\varphi_\h,\varphi_\g)$ being extensible.

\begin{defi}\label{ext}
Let $0\lon\h\stackrel{\iota}{\lon}\hat{\g}\stackrel{p}{\lon}\g\lon 0$ be a central extension of Lie algebras. A pair of derivations $(\varphi_\h,\varphi_\g)\in\Der(\h)\times\Der(\g)$ is said to be {\bf extensible} if there exists a derivation $\varphi_{\hat{\g}}\in \Der(\hat{\g})$
such that we have the following exact sequence of $\LD$ pair morphisms
\[\begin{CD}
0@>>>\mathfrak{h}@>\iota>>\hat{\mathfrak{g}}@>p>>{\mathfrak{g}}             @>>>0\\
@.    @V\varphi_\h VV   @V\varphi_{\hat{\mathfrak{g}}}VV  @V\varphi_\g VV    @.\\
0@>>>\mathfrak{h}@>\iota>>\hat{\mathfrak{g}}@>p>>{\mathfrak{g}}             @>>>0.
\end{CD}\]
Equivalently, $(\hat{\g},\varphi_{\hat{\g}})$ is a central extension of $(\g,\varphi_\g)$ by $(\h,\varphi_\h)$.
\end{defi}

     Let $s:\g\lon\hat{\g}$ be an arbitrary section of the central extension $0\lon\h\stackrel{\iota}{\lon}\hat{\g}\stackrel{p}{\lon}\g\lon 0$. Then any element of $\hat{\g}$ can be written uniquely as $s(x)+h$ for some $x\in\g$ and $h\in\h$. Define $\psi:\wedge^2\g\lon\h$ by
    \begin{equation}\label{eq:cen2cocycle}
      \psi(x,y)=[s(x),s(y)]_{\hat{\mathfrak{g}}}-s[x,y]_{\mathfrak{g}}.
    \end{equation}
   For any pair $(\varphi_\h,\varphi_\g)\in\Der(\h)\times\Der(\g)$, define a $2$-cochain $\obe\in C^2(\g;\h)$  by
\begin{eqnarray}\label{eq:extclass}
\obe(x,y)=\varphi_\h(\psi(x,y))-\psi(\varphi_\g(x),y)-\psi(x,\varphi_\g(y)),\,\quad \forall x,y\in\g.
\end{eqnarray}

\begin{pro}\label{pro:obe}
Let $0\lon\h\stackrel{\iota}{\lon}\hat{\g}\stackrel{p}{\lon}\g\lon 0$ be a central extension of Lie algebras. For any pair $(\varphi_\h,\varphi_\g)\in\Der(\h)\times\Der(\g)$, the $2$-cochain $\obe\in C^2(\g;\h)$ defined by \eqref{eq:extclass} is a $2$-cocycle of the Lie algebra $\g$ with the coefficient in the trivial representation $(\rho=0,\h)$. Moreover, the cohomology class $[\obe]\in \huaH^2(\g;\h)$ does not depend on the choice of sections.
\end{pro}

The cohomology class $[\obe]\in \huaH^2(\g;\h)$ is called the {\bf obstruction class} of $(\varphi_\h,\varphi_\g)$ being extensible.

\pf Let $s:\g\lon\hat{\g}$ be a section of the central extension of Lie algebras $0\lon\h\stackrel{\iota}{\lon}\hat{\g}\stackrel{p}{\lon}\g\lon 0$. The $\psi:\wedge^2\g\lon\h$ defined by \eqref{eq:cen2cocycle} is a $2$-cocycle of the Lie algebra $\g$ with the coefficient in the trivial representation $(\rho=0,\h)$, i.e.
$$\psi([x,y]_\g,z)+\psi([y,z]_\g,x)+\psi([z,x]_\g,y)=0.$$
Then we have
\begin{eqnarray*}
(\dM \obe)(x,y,z)&=&-\obe([x,y]_\g,z)-\obe([y,z]_\g,x)-\obe([z,x]_\g,y)\\
                 &=&-\varphi_\h(\psi([x,y]_\g,z))+\psi(\varphi_\g([x,y]_\g),z)+\psi([x,y]_\g,\varphi_\g(z))\\
                 &&-\varphi_\h(\psi([y,z]_\g,x))+\psi(\varphi_\g([y,z]_\g),x)+\psi([y,z]_\g,\varphi_\g(x))\\
                 &&-\varphi_\h(\psi([z,x]_\g,y))+\psi(\varphi_\g([z,x]_\g),y)+\psi([z,x]_\g,\varphi_\g(y))\\
                 &=&-\varphi_\h\big(\psi([x,y]_\g,z)+\psi([y,z]_\g,x)+\psi([z,x]_\g,y)\big)\\
                 &&+\psi([\varphi_\g(x),y]_\g,z)+\psi([y,z]_\g,\varphi_\g(x))+\psi([z,\varphi_\g(x)]_\g,y)\\
                 &&+\psi([x,\varphi_\g(y)]_\g,z)+\psi([\varphi_\g(y),z]_\g,x)+\psi([z,x]_\g,\varphi_\g(y))\\
                 &&+\psi([x,y]_\g,\varphi_\g(z))+\psi([y,\varphi_\g(z)]_\g,x)+\psi([\varphi_\g(z),x]_\g,y)\\
                 &=&0,
\end{eqnarray*}
which implies that $\obe$ is a $2$-cocycle of the Lie algebra  $\g$ with the coefficient in the trivial representation $(\rho=0,\h)$.

Let $s_1$ and $s_2$ be two different sections. Define $\phi:\mathfrak{g}\lon \mathfrak{h}$ by $\phi(x)=s_1(x)-s_2(x)$. Then we have $$\psi_{1}(x,y)=\psi_{2}(x,y)-\phi([x,y]_{\g}).$$ Moreover, we have
\begin{eqnarray*}
\obe^1(x,y)&=&\varphi_\h(\psi_1(x,y))-\psi_1(\varphi_\g(x),y)-\psi_1(x,\varphi_\g(y))\\
           &=&\varphi_\h(\psi_{2}(x,y))-\varphi_\h(\phi([x,y]_{\g}))-\psi_2(\varphi_\g(x),y)+\phi([\varphi_\g(x),y]_\g)\\&&-\psi_2(x,\varphi_\g(y))+\phi([x,\varphi_\g(y)]_\g)\\
           &=&\obe^2(x,y)-(\varphi_\h\circ\phi-\phi\circ\varphi_\g)([x,y]_\g)\\
           &=&\obe^2(x,y)+\dM (\varphi_\h\circ\phi-\phi\circ\varphi_\g)(x,y).
\end{eqnarray*}
Thus, we have $[\obe^1]=[\obe^2]\in \huaH^2(\g;\h)$. The proof is finished. \qed\vspace{3mm}

Now we are ready to give the main result in this section.

\begin{thm}\label{thm:derivationext}
Let $0\lon\h\stackrel{\iota}{\lon}\hat{\g}\stackrel{p}{\lon}\g\lon 0$ be a central extension of Lie algebras.   Then a pair $(\varphi_\h,\varphi_\g)\in\Der(\h)\times\Der(\g)$ is extensible if and only if the obstruction class $[\obe]\in \huaH^2(\g;\h)$ is trivial.
\end{thm}

\pf Let $s:\g\lon\hat{\g}$ be a section of the central extension of Lie algebras $0\lon\h\stackrel{\iota}{\lon}\hat{\g}\stackrel{p}{\lon}\g\lon 0$. Suppose that $(\varphi_\h,\varphi_\g)$ is extensible. Then there exists a derivation $\varphi_{\hat{\g}}\in \Der(\hat{\g})$ such that
we have the exact sequence of $\LD$ pair morphisms in Definition \ref{ext}. By $\varphi_\g\circ p=p\circ \varphi_{\hat{\g}}$, we obtain $\varphi_{\hat{\g}}(s(x))-s(\varphi_\g(x))\in\h$. Define $\lambda:\g\lon\h$ by $$\lambda(x)=\varphi_{\hat{\g}}(s(x))-s(\varphi_\g(x)).$$
Then we have
\begin{eqnarray*}
  \varphi_{\hat{\g}}(s(x)+h)&=&\varphi_{\hat{\g}}(s(x))+\varphi_\h(h)\\
                          &=&\varphi_{\hat{\g}}(s(x))-s(\varphi_\g(x))+s(\varphi_\g(x))+\varphi_\h(h)\\
                          &=&s(\varphi_\g(x))+\lambda(x)+\varphi_\h(h).
\end{eqnarray*}

Let $s(x)+h$ and $s(y)+l$ be any two elements of $\hat{\g}$. Since $\varphi_{\hat{\g}}$ is a derivation of $\hat{\g}$,
on one hand, we have
\begin{eqnarray*}
\varphi_{\hat{\g}}([s(x)+h,s(y)+l]_{\hat{\g}})&=&\varphi_{\hat{\g}}([s(x),s(y)]_{\hat{\g}})\\
                                              &=&\varphi_{\hat{\g}}(s[x,y]_\g+[s(x),s(y)]_{\hat{\g}}-s[x,y]_\g)\\
                                              &=&\varphi_{\hat{\g}}(s[x,y]_\g+\psi(x,y))\\
                                              &=&s(\varphi_\g([x,y]_\g))+\lambda([x,y]_\g)+\varphi_\h(\psi(x,y)).
\end{eqnarray*}
On the other hand, we have
\begin{eqnarray*}
&&[\varphi_{\hat{\g}}(s(x)+h),s(y)+l]_{\hat{\g}}+[s(x)+h,\varphi_{\hat{\g}}(s(y)+l)]_{\hat{\g}}\\
&=&[s(\varphi_\g(x))+\lambda(x)+\varphi_\h(h),s(y)+l]_{\hat{\g}}+[s(x)+h,s(\varphi_\g(y))+\lambda(y)+\varphi_\h(l)]_{\hat{\g}}\\
&=&[s(\varphi_\g(x)),s(y)]_{\hat{\g}}+[s(x),s(\varphi_\g(y))]_{\hat{\g}}\\
&=&s[\varphi_\g(x),y]_\g+[s(\varphi_\g(x)),s(y)]_{\hat{\g}}-s[\varphi_\g(x),y]_\g+s[x,\varphi_\g(y)]_\g
+[s(x),s(\varphi_\g(y))]_{\hat{\g}}-s[x,\varphi_\g(y)]_\g\\
&=&s[\varphi_\g(x),y]_\g+\psi(\varphi_\g(x),y)+s[x,\varphi_\g(y)]_\g+\psi(x,\varphi_\g(y)).
\end{eqnarray*}
Thus, we have   \begin{equation}\label{eq:ext}
  \varphi_\h(\psi(x,y))-\psi(\varphi_\g(x),y)-\psi(x,\varphi_\g(y))=-\lambda([x,y]_\g),
  \end{equation}
  which implies that
  $$
  \obe=\dM \lambda.
  $$
Therefore, the obstruction class is trivial.

Conversely, if the obstruction class is trivial, then there exists a $\lambda:\g\lon\h$ such that $\obe=\dM \lambda.$ For any element $s(x)+h\in\hat{\g},$   define $\varphi_{\hat{\g}}$ by
$$\varphi_{\hat{\g}}(s(x)+h)=s(\varphi_\g(x))+\lambda(x)+\varphi_\h(h).$$
By \eqref{eq:ext}, we obtain the exact sequence of $\LD$ pair morphisms in Definition \ref{ext}. Thus,  $(\varphi_\h,\varphi_\g)$ is extensible. The proof is finished. \qed \vspace{3mm}

Obviously, we have

\begin{cor}
Let $0\lon\h\stackrel{\iota}{\lon}\hat{\g}\stackrel{p}{\lon}\g\lon 0$ be a central extension of Lie algebras. If $ \huaH^2(\g;\h)=0$, then any pair $(\varphi_\h,\varphi_\g)$ in $\Der(\h)\times\Der(\g)$ is extensible.
\end{cor}

At the end of this section, we give the condition on a pair of derivations $(\varphi_\h,\varphi_\g)\in\Der(\h)\times\Der(\g)$ such that it is extensible in every central extension of $\g$ by $\h$.
 By Proposition \ref{pro:obe}, we can define a linear map $\Theta:\Der(\h)\times\Der(\g)\lon \gl(\huaH^2(\g;\h))$ by
\begin{eqnarray*}
\Theta(\varphi_\h,\varphi_\g)([\psi])=[\varphi_\h\circ \psi-\psi(\varphi_\g\otimes{\Id})-\psi({\Id}\otimes\varphi_\g)].
\end{eqnarray*}

\begin{thm}
Let $\h$ be an abelian Lie algebra and $\g$ a Lie algebra. A pair of derivations $(\varphi_\h,\varphi_\g)\in\Der(\h)\times\Der(\g)$ is extensible in every central extension of $\g$ by $\h$ if and only if $\Theta(\varphi_\h,\varphi_\g)=0$.
\end{thm}

\pf We suppose $\Theta(\varphi_\h,\varphi_\g)=0$. For any central extension $0\lon\h\stackrel{\iota}{\lon}\hat{\g}\stackrel{p}{\lon}\g\lon 0$, we choose a section $s:\g\lon\hat{\g}$. Then $\psi:\wedge^2\g\lon\h$ defined by $\psi(x,y)=[s(x),s(y)]_{\hat{\mathfrak{g}}}-s[x,y]_{\mathfrak{g}}$ is a $2$-cocycle. Moreover, we obtain that
\begin{eqnarray*}
[\obe]=[\varphi_\h\circ \psi-\psi(\varphi_\g\otimes{\Id})-\psi({\Id}\otimes\varphi_\g)]=\Theta(\varphi_\h,\varphi_\g)([\psi])=0.
\end{eqnarray*}
By Theorem \ref{thm:derivationext},
  $(\varphi_\h,\varphi_\g)$ is extensible in the above central extension.

Conversely, for any element $[\psi]\in \huaH^2(\g;\h)$, there exists a central extension $0\lon\h\stackrel{\iota}{\lon}\g\oplus \h\stackrel{p}{\lon}\g\lon 0$, where the bracket on $\g\oplus\h$ is defined by
$$
[x+g,y+h]=[x,y]_\g+\psi(x,y).
$$ Since $(\varphi_\h,\varphi_\g)$ is extensible in every central extension of $\g$ by $\h$, by Theorem \ref{thm:derivationext}, we have
$$\Theta(\varphi_\h,\varphi_\g)([\psi])= [\varphi_\h\circ \psi-\psi(\varphi_\g\otimes{\Id})-\psi({\Id}\otimes\varphi_\g)]=0.$$
Therefore, we have $\Theta(\varphi_\h,\varphi_\g)=0$. The proof is finished. \qed

\section{Classification of skeletal $\Ld$ pairs}\label{2-Lie}

In this section, we call a Lie 2-algebra with a derivation of degree 0 a $\Ld$ pair and show that the third cohomology group  $\huaH^3_{\LD}(\g;V)$ classifies skeletal $\Ld$ pairs. See \cite{Baze} for more details about Lie 2-algebras.

\begin{defi}\label{defi:Lie2}
A Lie $2$-algebra $\mathcal{V}$ consists of the following data:
\begin{itemize}
\item a complex of vector spaces $V_{1}\stackrel{l_1}{\longrightarrow} V_{0}$,
    \item bilinear maps $l_{2}:V_{i}\times V_{j}\longrightarrow V_{i+j}$, where $0 \le i + j \le 1,$
       \item a skew-symmetric trilinear map $l_{3}:V_{0}\times V_{0}\times V_{0}\longrightarrow V_{1}$,
      such that for all $w,x,y,z\in V_{0}$ and $m,n\in V_{1}$, the following equalities are satisfied:
        \item[\rm(a)] $l_{2}(x,y)=-l_{2}(y,x)$,\quad $l_{2}(x,m)=-l_{2}(m,x)$,
        \item[\rm(b)]$l_1 l_{2}(x,m)=l_{2}(x,l_1 m)$,\quad $l_{2}(l_1 m,n)=l_{2}(m,l_1 n)$,
        \item[\rm(c)]$l_1 l_{3}(x,y,z)=l_{2}(x,l_{2}(y,z))+l_{2}(y,l_{2}(z,x))+l_{2}(z,l_{2}(x,y))$,
        \item[\rm(d)]$l_{3}(x,y,l_1 m)=l_{2}(x,l_{2}(y,m))+l_{2}(y,l_{2}(m,x))+l_{2}(m,l_{2}(x,y))$,
        \item[\rm(e)]
     $
        l_{3}(l_{2}(w,x),y,z)+l_{2}(l_{3}(w,x,z),y)
        +l_{3}(w,l_{2}(x,z),y)\\
        +l_{3}(l_{2}(w,z),x,y)
        =l_{2}(l_{3}(w,x,y),z)+l_{3}(l_{2}(w,y),x,z)\\+l_{3}(w,l_{2}(x,y),z)
        +l_{2}(w,l_{3}(x,y,z))+l_{2}(l_{3}(w,y,z),x)+l_{3}(w,l_{2}(y,z),x).
     $
\end{itemize}
\end{defi}
We usually denote a Lie 2-algebra by $(V_{1},V_{0},l_1,l_{2},l_{3})$ or simply by $\mathcal{V}$. A Lie 2-algebra is called  {\bf skeletal} if $l_1=0$. There is a one-to-one correspondence between skeletal Lie 2-algebras and triples $(\g,(\rho,V),l_3)$, where $\g$ is a Lie algebra, $(\rho,V)$ is a representation of $\g$, and $l_3$ is a 3-cocycle on $\g$ with the coefficient in $V$. More precisely,  for a skeletal Lie $2$-algebra $\mathcal{V}=(V_{1},V_{0},l_1=0,l_{2},l_{3})$, it is straightforward to see that $(V_{0},l_{2})$ is a Lie algebra. Define $\rho$ from $V_{0}$ to $\gl(V_{1})$ by
\begin{eqnarray*}
\rho(x)(u)=l_2(x,u).
\end{eqnarray*}
Then, $(\rho,V_{1})$ is a representation of the Lie algebra $(V_0,l_{2})$ and $l_3$ is a 3-cocycle on $V_{0}$ with the coefficient in $V_{1}$.

\begin{defi}
Let $\mathcal{V}$ and $\mathcal{V'}$ be Lie $2$-algebras. A {\bf morphism} $f$ from $\mathcal{V}$ to $\mathcal{V'}$ consists of:
\begin{itemize}
  \item a chain map $f:\mathcal{V}\lon \mathcal{V'}$, which consists of linear maps $f_{0}:V_{0}\lon V_{0}'$ and $f_{1}:V_{1}\lon V_{1}'$ satisfying $$f_{0}\circ l_1=l_1'\circ f_{1},$$
  \item a skew-symmetric bilinear map $f_{2}:V_{0}\times V_{0}\lon V_{1}'$ such that for all $x,y,z\in V_{0}$ and $m\in V_{1}$,
  \end{itemize}
  the following equalities hold:
  \begin{itemize}
  \item $l_1' f_{2}(x,y)=f_{0}(l_{2}(x,y))-l_{2}'(f_{0}(x),f_{0}(y))$,
  \item $f_{2}(x,l_1 m)=f_{1}(l_{2}(x,m))-l_{2}'(f_{0}(x),f_{1}(m))$,
  \item
$
   l_{2}'(f_{0}(x),f_{2}(y,z))+l_{2}'(f_{0}(y),f_{2}(z,x))
   +l_{2}'(f_{0}(z),f_{2}(x,y))+l_{3}'(f_{0}(x),f_{0}(y),f_{0}(z))\\
    =f_{2}(l_{2}(x,y),z)+f_{2}(l_{2}(y,z),x)+f_{2}(l_{2}(z,x),y)+f_{1}(l_{3}(x,y,z)).
$
\end{itemize}
\end{defi}
A morphism is called an isomorphism if $f_0$ and $f_1$ are invertible.

\begin{defi}\label{defi:derivation}
A derivation of degree $0$ of a Lie $2$-algebra $\mathcal{V}$ is a triple $(X_0,X_1,l_X)$, in which $X=(X_0,X_1)\in \End(V_0)\oplus\End(V_1)$ and $l_X:V_0\wedge V_0\lon V_1$ is a linear
map, such that for all $x,y,z\in V_0,m\in V_1$, the following equalities hold:
\begin{itemize}
\item[\rm(a)]
$X_0\circ l_1=l_1\circ X_1,$
\item[\rm(b)]
$l_1 l_X(x,y)=X_0(l_2(x,y))-l_2(X_0x,y)-l_2(x,X_0y),$
\item[\rm(c)]
$l_X(x,l_1 m)=X_1(l_2(x,m))-l_2(X_0x,m)-l_2(x,X_1m),$
\item[\rm(d)]
$X_1l_3(x,y,z)=l_X(x,l_2(y,z))+l_2(x,l_X(y,z))+l_3(X_0x,y,z)+c.p.(x,y,z).$
\end{itemize}
\end{defi}

See \cite{Doubek-Lada,Lang-Liu-Sheng} for more details about derivations of Lie 2-algebras and $L_\infty$-algebras. We denote a Lie 2-algebra with a derivation of degree 0 by $(\huaV;(X_0,X_1,l_X))$ and call it a $\Ld$ pair. In particular, a skeletal Lie 2-algebra with a derivation of degree 0 will be called a skeletal $\Ld$ pair.

\begin{defi}\label{defi:derivationiso}
Let $(\huaV;(X_0,X_1,l_X))$ and $(\huaV';(X_0',X_1',l_X'))$ be $\Ld$ pairs. An isomorphism $\tau$ from $(\huaV;(X_0,X_1,l_X))$ to $(\huaV';(X_0',X_1',l_X'))$ consists of linear maps $f_0:V_0\lon V_0',f_1:V_1\lon V_1',$ $f_2:V_0\wedge V_0\lon V_1'$ and $\huaB:V_0\lon V_1'$ such that $(f_0,f_1,f_2)$ is a Lie $2$-algebra isomorphism from $\mathcal{V}$ to $\mathcal{V'}$ and the following equalities hold for all $x,y\in V_0,m\in V_1$:
\begin{itemize}
\item[\rm(a)]
$X_0'(f_0(x))-f_0(X_0(x))=l_1'(\huaB(x)),$
\item[\rm(b)]
$X_1'(f_1(m))-f_1(X_1(m))=\huaB(l_1(m)),$
\item[\rm(c)]
$f_1(l_X(x,y))+f_2(X_0x,y)+f_2(x,X_0y)-X_1'(f_2(x,y))-l_X'(f_0(x),f_0(y))=l_2'(\huaB(x),f_0(y))+l_2'(f_0(x),\huaB(y))-\huaB(l_2(x,y)).$
\end{itemize}
\end{defi}
\emptycomment{
\begin{rmk}
Let $((\g,[\cdot,\cdot]_\g);\varphi_\g)$ and $((\g',[\cdot,\cdot]_{\g'});\varphi_{\g'})$ be $\LD$ pairs. We say $((\g,[\cdot,\cdot]_\g);\varphi_\g)$ and $((\g',[\cdot,\cdot]_{\g'});\varphi_{\g'})$ are isomorphism if and only if there is a Lie algebra isomorphism $f:\g\lon\g'$ such that $f\circ\varphi_\g=\varphi_{\g'}\circ f$, which is equivalent to
\begin{eqnarray*}
\frac{d}{dt}\mid_{t=0}\Big(f\circ \exp(t\varphi_\g)-\exp(t\varphi_{\g'})\circ f\Big)=0.
\end{eqnarray*}
Moreover, let $(\huaV;(X_0,X_1,l_X))$ and $(\huaV';(X_0',X_1',l_X'))$ be $\Ld$ pairs. There is a Lie $2$-algebra isomorphism $f=(f_0,f_1,f_2):\huaV\lon\huaV'$. By integration of the derivation of Lie $2$-algebra, we have Lie $2$-algebra automorphisms
$\exp(t,(X_0,X_1,l_X)):\huaV\lon\huaV$ and $\exp(t,(X_0',X_1',l_X')):\huaV'\lon\huaV'$.
\end{rmk}}

Let $\mathcal{V}=(V_{1},V_{0},l_{1}=0,l_{2},l_{3})$ be a skeletal Lie $2$-algebra and $(X_0,X_1,l_X)$ be a derivation of degree 0 of $\mathcal{V}$.
By
condition (b) in Definition \ref{defi:derivation}, we deduce that $X_0$ is a derivation of the Lie algebra $(V_{0},l_{2})$. Condition (c) in Definition \ref{defi:derivation} implies that $(\rho,V_1,X_1)$ is a representation of the $\LD$ pair $(V_{0},X_0)$. Condition (d) in Definition \ref{defi:derivation} implies that $(l_3,-l_X)\in \huaZ^3_{\LD}(V_0;V_1)$.

Conversely, if $(\rho,V,\varphi_V)$ is a representation of the $\LD$ pair $(\g,\varphi_\g)$ and $(\theta_3,\theta_2)\in \huaZ^3_{\LD}(\g;V)$, then we can deduce that $(\varphi_\g,\varphi_V,-\theta_2)$ is a derivation of the skeletal Lie 2-algebra $(V\stackrel{0}{\rightarrow}\g,l_2,l_3=\theta_3)$, where $l_2$ is defined by
\begin{eqnarray*}
l_2(x,y)=[x,y]_\g, \quad l_2(x,u)=-l_2(u,x)=\rho(x)(u),\quad \forall x,y\in \g, u\in V.
\end{eqnarray*}
Summarize the above discussion, we have

\begin{pro}\label{correspondence}
There is a one to one correspondence between skeletal $\Ld$ pairs and triples $\big((\g,\varphi_\g),(\rho,V,\varphi_V),(\theta_3,\theta_2)\big)$, where $(\g,\varphi_\g)$ is a $\LD$ pair , $(\rho,V,\varphi_V)$ is a representation of $(\g,\varphi_\g)$, and $(\theta_3,\theta_2)$ is a $3$-cocycle on $(\g,\varphi_\g)$ with the coefficient in $(\rho,V,\varphi_V)$.
\end{pro}

In the sequel, we give the equivalence relation between triples $\big((\g,\varphi_\g),(\rho,V,\varphi_V),(\theta_3,\theta_2)\big)$ and show that there is a one-to-one correspondence between equivalence classes of such triples and isomorphism classes of skeletal $\Ld$ pairs.

\begin{defi}\label{equivalent}
 Let $\big((\g,\varphi_\g),(\rho,V,\varphi_V),(\theta_3,\theta_2)\big)$ and $\big((\g',\varphi_{\g'}),(\rho',V',\varphi_{V'}),(\theta_3',\theta_2')\big)$ be triples as described in Proposition \ref{correspondence}. They are said to be equivalent if there exist Lie algebra isomorphism $\alpha:\g\lon\g'$, linear isomorphism $\beta:V\lon V'$ and two linear maps $\gamma:\g\wedge\g\lon V',\,\,\eta:\g\lon V'$ such that the following equalities hold for all $x,y,z\in\g,u\in V$:
 \begin{itemize}
\item[\rm(a)]
$\varphi_{\g'}(\alpha(x))=\alpha(\varphi_\g(x)),$
\item[\rm(b)]
$\varphi_{V'}(\beta(u))=\beta(\varphi_V(u)),$
\item[\rm(c)]
$\beta(\rho(x)(u))=\rho'(\alpha(x))(\beta(u))$,
\item[\rm(d)]
$\rho'(\alpha(x))(\gamma(y,z))+\rho'(\alpha(y))(\gamma(z,x))
   +\rho'(\alpha(z))(\gamma(x,y))+\theta_3'(\alpha(x),\alpha(y),\alpha(z))\\
    =\gamma([x,y]_\g,z)+\gamma([y,z]_\g,x)+\gamma([z,x]_\g,y)+\beta(\theta_3(x,y,z)),$
\item[\rm(e)]
$-\beta(\theta_2(x,y))+\gamma(\varphi_\g (x),y)+\gamma(x,\varphi_\g (y))-\varphi_{V'}(\gamma(x,y))+\theta_2'(\alpha(x),\alpha(y))\\=
-\rho'(\alpha(y))(\eta(x))+\rho'(\alpha(x))(\eta(y))-\eta([x,y]_\g).$
\end{itemize}
\end{defi}

\begin{thm}
There is a one-to-one correspondence between isomorphism classes of skeletal $\Ld$ pairs and equivalence classes of triples $\big((\g,\varphi_\g),(\rho,V,\varphi_V),(\theta_3,\theta_2)\big)$, where $(\g,\varphi_\g)$ is a $\LD$ pair, $(\rho,V,\varphi_V)$ is a representation of $(\g,\varphi_\g)$, and $(\theta_3,\theta_2)$ is a $3$-cocycle on $(\g,\varphi_\g)$ with the coefficient in $(\rho,V,\varphi_V)$.
\end{thm}
\pf Let $\big((\g,\varphi_\g),(\rho,V,\varphi_V),(\theta_3,\theta_2)\big)$ and $\big((\g',\varphi_{\g'}),(\rho',V',\varphi_{V'}),(\theta_3',\theta_2')\big)$ be equivalent triples. By Proposition \ref{correspondence}, we have two skeletal Lie 2-algebras, given by
$$\huaV=(V\stackrel{0}{\rightarrow}\g,l_2,l_3=\theta_3),\,\,\,\,\huaV'=(V'\stackrel{0}{\rightarrow}\g',l_2',l_3'=\theta_3').$$
Moreover, $(X_0=\varphi_\g,X_1=\varphi_V,l_X=-\theta_2)$ and $(X_0'=\varphi_{\g'},X_1'=\varphi_{V'},l_X'=-\theta_2')$ are  degree 0 derivations of $\huaV$ and $\huaV'$ respectively.

We define $f=(f_0=\alpha,f_1=\beta,f_2=\gamma)$, and $\huaB=\eta$. By condition (c) and condition (d) in Definition \ref{equivalent} and the fact that $\alpha$ is a Lie algebra isomorphism,   $f$ is a Lie $2$-algebra isomorphism from $\mathcal{V}$ to $\mathcal{V'}$. Moreover, by conditions (a), (b) and   (e) in Definition \ref{equivalent}, we deduce that $(\huaV;(X_0,X_1,l_X))$ is isomorphic to $(\huaV';(X_0',X_1',l_X'))$.

The converse part can be proved similarly and we omit details. \qed

Tang Rong and Yunhe Sheng\\
Department of Mathematics, Jilin University,
 Changchun 130012,  China\\
 Email: shengyh@jlu.edu.cn, \quad tangrong16@mails.jlu.edu.cn\\
Yael Fr\'egier\\
LML, Artois University,
 Lens 62307, France\\
 Email: yael.fregier@gmail.com


\begin{thebibliography}{999}


\bibitem{Ay-Ki-Tr}
V. Ayala, E. Kizil and I. de~Azevedo~Tribuzy,
\newblock On an algorithm for finding derivations of {L}ie algebras.
\newblock {\em Proyecciones} 31 (2012), 81-90.

\bibitem{Ay-Ti}
V. Ayala and J. Tirao,
\newblock Linear control systems on {L}ie groups and controllability.
\newblock In {\em Differential geometry and control ({B}oulder, {CO}, 1997)},
  volume~64 of {\em Proc. Sympos. Pure Math.} pages 47-64. Amer. Math. Soc.
  Providence, RI, 1999.


\bibitem{Barr-Beck} M. Barr and J. Beck, Homology and standard constructions, 1969 \emph{Sem. on Triples and Categorical Homology Theory (ETH, Zurich, 1966/67)} pp. 245-335 Springer, Berlin

\bibitem{Barr} M. Barr, Cartan-Eilenberg cohomology and triples, \emph{J. Pure Appl. Algebra} 112 (1996), no. 3, 219-238.

\bibitem{Ba-Vi}
I.~A.~Batalin and G.~A.~Vilkovisky,
\newblock Gauge algebra and quantization.
\newblock {\em Phys. Lett. B} 102 (1981), 27-31.


\bibitem{Baze} J. C. Baez and A. S. Crans, Higher-dimensional algebra. VI. Lie 2-algebras, \emph{Theory Appl. Categ.} 12 (2004), 492-538.

\bibitem{Bardakov-Singh}
V. G. Bardakov and M. Singh, Extensions and automorphisms of Lie algebras, \emph{J. Algebra Appl.} 16 (2017), 15 pp.




\bibitem{Ch-Ei}
C. Chevalley and S. Eilenberg,
\newblock Cohomology theory of {L}ie groups and {L}ie algebras.
\newblock {\em Trans. Amer. Math. Soc.} 63 (1948), 85-124.




\bibitem{Doubek-Lada}
M. Doubek and T. Lada, Homotopy derivations, \emph{J. Homotopy Relat. Struct.} 11 (2016), no. 3, 599-630.



\bibitem{Doubek-Markl-Zima}
M. Doubek, M. Markl and P. Zima, Deformation theory (lecture notes), \emph{ Arch. Math. (Brno)} 43 (2007), no. 5, 333-371.

\bibitem{Fregier-Markl-Yau}
Y. Fr\'egier, M. Markl and D. Yau, The $L_\infty$-deformation complex of diagrams of algebras, \emph{New York J. Math. } 15 (2009), 353-392.

\bibitem{Fregier}
Y. Fr\'egier, A new cohomology theory associated to deformations of Lie algebra morphisms, \emph{Lett. Math. Phys.} 70 (2004), no. 2, 97-107.


\bibitem{Fregier-Zambon1}
Y. Fr\'egier and M. Zambon, Simultaneous deformations and Poisson geometry, \emph{Compos. Math.} 151 (2015), no. 9, 1763-1790.

\bibitem{Fregier-Zambon2}
Y. Fr\'egier and M. Zambon, Simultaneous deformations of algebras and morphisms via derived brackets, \emph{J. Pure Appl. Algebra} 219 (2015), no. 12, 5344-5362.

\bibitem{Gerstenhaber1}
M. Gerstenhaber, The cohomology structure of an associative ring, \emph{Ann. of Math. (2)} 78 (1963), 267-288.

\bibitem{Gerstenhaber2}
M. Gerstenhaber, On the deformation of rings and algebras, \emph{Ann. of Math. (2) } 79 (1964), 59-103.

\bibitem{Gerstenhaber-Schack}
M. Gerstenhaber and S. D. Schack, On the deformation of algebra morphisms and diagrams, \emph{Trans. Amer. Math. Soc.} 279 (1983), no. 1, 1-50.


\bibitem{Gi-Ka}
V. Ginzburg and M. Kapranov,
\newblock Koszul duality for operads.
\newblock {\em Duke Math. J.} 76 (1994), 203-272.

\bibitem{Har}
D.~K. Harrison,
\newblock Commutative algebras and cohomology,
\newblock {\em Trans. Amer. Math. Soc.}, 104 (1962), 191-204.

\bibitem{Hor}
G.~Hochschild,
\newblock On the cohomology groups of an associative algebra.
\newblock {\em Ann. of Math. (2)}, 46 (1945), 58-67.




\bibitem{Lang-Liu-Sheng} H. Lang, Z. Liu and Y. Sheng, Integration of derivations for Lie 2-algebras, \emph{Transform. Groups } 21 (2016), no. 1, 129-152.

\bibitem{Loday}
L. J. Loday, On the operad of associative algebras with derivation, \emph{Georgian Math. J.} 17 (2010), no. 2, 347-372.

\bibitem{Markl}
M. Markl, Intrinsic brackets and the
$L_{\infty}$-deformation theory of bialgebras, \emph{J. Homotopy Relat. Struct.} 5 (2010), no. 1, 177-212.

\bibitem{Mandal}
A. Mandal, Deformation of Leibniz algebra morphisms, \emph{Homology Homotopy Appl.} 9 (2007), no. 1, 439-450.

\bibitem{Nijenhuis-Richardson}
A. Nijenhuis and R. Richardson, Cohomology and deformations in graded Lie algebras, \emph{Bull. Amer. Math. Soc.} 72 (1966), 1-29.



\bibitem{Sta}
J. Stasheff,
\newblock Homotopy associativity of {$H$}-spaces. {I}, {II}.
\newblock {\em Trans. Amer. Math. Soc.} 108 (1963), 275-292.

\bibitem{Vai}
A.~Yu. Vaintrob,
\newblock Lie algebroids and homological vector fields,
\newblock {\em Uspekhi Mat. Nauk}, 52 (1997), 161-162.

\bibitem{Vor}
Th.~Th. Voronov,
\newblock Higher derived brackets for arbitrary derivations,
\newblock In {\em Travaux math{\'e}matiques. {F}asc. {XVI}}, volume~16 of {\em
  Trav. Math.}, pages 163-186. Univ. Luxemb., Luxembourg, 2005.



\bibitem{Yau}
D. Yau, Deformations of coalgebra morphisms, \emph{J. Algebra} 307 (2007), no. 1, 106-115.
\end{thebibliography}
\end{document}